\font\title=cmr10 scaled 1200

\newcount\notenumber

\def\note{\advance\notenumber by 1
\footnote{$^{(\the\notenumber)}$}}

\def\note{\advance\notenumber by 1
\footnote{$^{(\the\notenumber)}$}}

\font\tenmsam=msam10
\font\tenmsb=msbm10
\font\sevenmsb=msbm7
\font\fivemsb=msbm5
\newfam\msbfam
\textfont\msbfam=\tenmsb
\scriptfont\msbfam=\sevenmsb
\scriptscriptfont\msbfam=\fivemsb
\def\quadratino{{\tenmsam\char003}}
\def\cvd{{\unskip\nobreak\penalty50
\hskip1em\hbox{ }\hskip2em\hbox{}\nobreak\hfil\quadratino
\parfillskip=0pt \finalhyphendemerits=0 \par}
\vskip8pt plus4pt minus2pt}

\font\tengoth=eufm10
\font\sevengoth=eufm7
\font\fivegoth=eufm5

\newfam\gothfam \scriptscriptfont\gothfam=\fivegoth
\textfont\gothfam=\tengoth \scriptfont\gothfam=\sevengoth

\def\N{{\bf N}}

\def\Z{{\bf Z}}

\def\C{{\bf C}}
\def\Q{{\bf Q}}

\def\Pr{{\bf P}}
\def\Gm{{\bf G}_m}

\def\a{{\bf a}}
\def\b{{\bf b}}

\def\O{{\cal O}}
\def\c{{\cal C}}

\def\cc{{\tilde{{\cal C}}}}
\def\dd{{\tilde{{\cal D}}}}

\def\A{{\bf A}}

\def\uno{{u_1^\prime\over u_1}}
\def\due{{u_2^\prime\over u_2}}

\centerline{\title
Some cases of Vojta's Conjecture on integral points over function fields}

\bigskip

\centerline{{\sl Pietro Corvaja\quad \&\quad Umberto Zannier}}

\bigskip

$${}$$

\noindent{\bf Abstract}. In the present paper we solve in particular  the  function field
version of  a special case of Vojta's conjecture for integral points, namely for the variety
obtained by removing a conic and two lines from the projective plane. This will follow  from  a
bound for the degree of a  curve on such a surface in terms of its Euler characteristic.

This case is special, but significant, because it lies ``at the boundary", in the sense that it
represents the simplest case of the conjecture which is still open. Also, it was studied in the
context of Nevanlinna Theory by M. Green already in the seventies.

Our general  results immediately imply the degeneracy of solutions of
Fermat's type equations $z^d=P(x^m,y^n)$ for all $d\ge 2$ and large enough $m,n$, also in the
case of non-constant coefficients.  Such equations fall apparently out of all known treatments.

The methods used here refer to derivations, as is usual in
function fields, but contain  fundamental new points. One of the tools concerns an estimation for
the $\gcd (1-u,1-v)$ for $S$-units $u,v$; this had been developed also in the arithmetic case,
 but for function fields  we may obtain a much more uniform quantitative version.

In an Appendix we shall finally  point out some other implications of the methods  to the problem
of torsion-points on curves and related known questions.
\bigskip

\noindent{\bf \S1. Introduction and main results}.
\smallskip

A celebrated conjecture in
diophantine geometry, first proposed by Vojta (see [V, Conjecture 3.4.3] and [V, Prop. 4.1.2]),
reads as follows:

{\it Let $X$ be a smooth affine variety defined over a number field $k$,
$\tilde X$ be a smooth projective variety containing $X$ as an open subset,
$D=\tilde{X}\setminus X$ the divisor at infinity and $K$ a canonical
divisor of $\tilde{X}$. Suppose that $D$ is a normal crossing divisor.
Then if $D+K$ has maximal Kodaira dimension, for every ring of $S$-integers
$\O_S\subset k$, the set of $S$-integral points $X(\O_S)$ is not Zariski-dense.}
\smallskip

It is known that the Kodaira dimension of $D+K$ is in fact independent of the
smooth compactification $\tilde{X}$ of $X$, provided that $D$ has
normal crossings (see for instance [KMK]). Following a
frequent notation, we will call log Kodaira dimension of $X$ the Kodaira dimension
of $D+K$.
\smallskip

A complex analytic analogue  of Vojta's Conjecture asks for the degeneracy
of entire curves on affine varieties with maximal log Kodaira
dimension; more precisely, one expects that, for every holomorphic map
$f:\C\rightarrow X(\C)$ to such a variety,
the image $f(\C)$ be contained in a proper closed algebraic subvariety.
\smallskip

Finally an (algebraic)  function field  analogue of Vojta's Conjecture predicts that,
given a smooth curve $\cc$ and a finite subset $S\subset\cc$,
there should exist a bound for the degree (in a suitable projective embedding)
for the images of non-constant morphisms  $\cc\setminus S\rightarrow X$, where $X$ is again an
algebraic variety with maximal log Kodaira dimension.  \smallskip

The particular case where $\tilde{X}$ is the projective plane has been widely studied,
in the arithmetic, analytic and algebraic setting.
The condition on the Kodaira dimension of $D+K$ is equivalent to the inequality $\deg(D)\geq 4$.
In all settings, it turns out that the case where $D$ has at least four components is much easier:
in this case, and only in this case, $X$ admits a finite map to  a closed subvariety of
the torus $\Gm^3$,  that   is not  the translate of a subtorus.
Then the degeneracy of integral points is equivalent to the three-variable $S$-unit Theorem, which
follows easily from the   Subspace Theorem in diophantine approximation; the degeneracy of entire
curves  is a consequence of a famous theorem of Emile Borel (see [Bo]);  finally, the function
field analogue is obtained as an application of the Stothers-Mason $abc$ Theorem (see [BrM]).
\smallskip

Hence the first crucial open case occurs with $\deg(D)=4$ and $D$ consisting of three
irreducible components, i.e. $D$ is the sum of a conic and two lines.  In the arithmetic case,
first considered by Beukers [Be, p. 116], the solution  seems to be, at present, beyond hope  (but
some particular cases have been settled in [CZ1]). An open question, which would have a positive
answer under Vojta's conjecture for the complement of a conic and two lines, is the  finiteness of
the perfect squares of the form $2^a+3^b+1$ for positive integers $a,b$. 

The Nevanlinna-theoretic analogue has been investigated since the seventies: Mark Green attacked a
particular case,  which amounts to the characterisation of squares,  in the ring of entire
functions, of the form $\exp(\varphi)+\exp(\psi)+1$  for entire functions $\varphi,\psi$; he
characterises  such functions under technical hypothesis on the growth of $\varphi,\psi$.
Recently, Noguchi,  Winkelmann and Yamanoi [NWY] solved completely the  problem, even in much
greater generality, by classifying in particular the entire curves on the complement of a conic and
two lines. \smallskip

The present paper solves the (algebraic) function field analogue, which to us seems to be
nontrivial already in the case of rational curves.
\smallskip

We now describe precisely our main result. We let $\kappa$ be an algebraically closed field of
characteristic zero; all algebraic varieties will be defined over $\kappa$. We let
$X=\Pr_2\setminus D$, where $D$ is a quartic consisting of the  union of a smooth conic and two
lines in general position (i.e. $D$ has five singular points).  Let $\c\subset X$ be an affine
curve on $X$;  its normalization is of the form $\cc\setminus S$, for a (unique) smooth complete
curve $\cc$ and a finite subset $S\subset\cc$. We define the {\it Euler  characteristic} of $\c$
to be the Euler characteristic of $\cc\setminus S$, i.e. $$ \chi(\c)=\chi_S(\cc)=2g(\cc)-2+\sharp
S. $$ This definition is consistent with Def. 1.2 of [KMK], where it is said that
a curve $\bar{\c}$ meets $k$ times a divisor $D$ if the normalization of
$\bar{\c}\setminus (D\cap\bar{\c})$ has $k$ points at infinity.

Finally, we let $\deg(\c)$ be the projective degree of the projective closure
$\bar{\c}$ of $\c$ in $\Pr_2$.
Our first result is
\medskip

\noindent{\bf Theorem 1.1}. {\it Let, as before, $X=\Pr_2\setminus D$, where $D$ is a quartic consisting of the union of a smooth conic and two lines in general position. Let
$\c\subset X$ be an affine curve. Then
$$
\deg (\c)\leq 2^{15}\cdot 35\cdot \max\{1,\chi(\c)\}. \eqno(1.1)
$$
}
\bigskip

Theorem 1.1 will be proved after reducing it (via Lemma 3.3) to the ``diophantine equation"
$$
y^2=u_1^2+\lambda u_1+u_2+1,\eqno(1.2)
$$
to be solved in rational functions $(y,u_1,u_2)\in\kappa(\cc)$
where $y$ is regular on $\cc\setminus S$ and $u_1,u_2$ are $S$-units, i.e. regular and never
vanishing functions on  $\cc\setminus S$ (see Theorem 3.4).
A natural generalization, which will be obtained here by the same method,
concerns more general bounds for the number of multiple zeros of rational functions of the form
$A(u_1,u_2)$ for polynomials $A(X,Y)\in\kappa[X,Y]$.  We shall prove the following generalization
of Theorem 3.4, where the symbol $H(\cdot)$ is used to denote the {\it height} of a rational
function on $\cc$, i.e. its degree as a morphism $\cc\rightarrow\Pr_1$.  \medskip

\noindent{\bf Theorem 1.2}. {\it Let $\cc$ be a smooth complete curve,
$S\subset\cc$ a finite set of points. Let $A(X,Y)\in\kappa[X,Y]$ be a
polynomial without repeated factors, $\epsilon>0$ a positive real number.  
There exists a number $C_1=C_1(A,\epsilon)$ such that 
for all pairs $(u_1,u_2)\in(\O_S^*)^2$ with $\max\{H(u_1),H(u_2)\}\geq
C_1\cdot\max\{1,\chi(\cc\setminus S)\}$, either the rational function $A(u_1,u_2)$ satisfies
$$
\sum_{v\in\cc\setminus S} \max\{0,v(A(u_1,u_2)-1)\}<\epsilon\max\{H(u_1),H(u_2)\}, \eqno(1.3)
$$
 or the $S$-units  $u_1,u_2$ verify a multiplicative dependence relation of the form
$$
u_1^r\cdot  u_2^s\in\kappa\eqno(1.4)
$$
for a suitable pair of integers $(r,s)\in\Z^2\setminus\{0\}$.
}
 
\medskip

Plainly, Theorem 1.2 applies in particular to the equation
$$
y^d=A(u_1,u_2),\qquad d\ge 2,\eqno(1.5)
$$
to be solved in $S$-units $u_1,u_2\in\O_S^*$ and $S$-integers $y\in\O_S$. As a result, we obtain
a bound for the height of solutions to (1.5), outside some possible infinite families satisfying
(1.4). Note that (1.2) is a particular case of (1.5); hence Theorem 1.2 generalizes Theorem 1.1.
However, in the proof of Theorem 1.2 we did not calculate all the numerical constants involved.

\smallskip

An immediate corollary of Theorem 1.2 (in the above quantification), which we do not state for
simplicity, is the degeneracy of the solutions for a Fermat-type equation
$$
z^d=A(x^n,y^m),
$$
for any $d\ge 2$ and $n,m$ large enough in  terms of $\deg A$. To get a proof, it suffices to
note that $x^n,y^m$ are $S$-units for a set $S$ which is small with respect to the involved
degrees.\medskip

Any bound for the height of the solutions to such equations (1.5), in turn,
corresponds to a bound for the degree of an affine curve, in terms of
its Euler characteristic, on affine varieties obtained as ramified cyclic covers of $\Gm^2$.
Note that the complement of a conic and two lines in $\Pr_2$ admits a finite map
of degree two to $\Gm^2$, i.e. $X$ is a double ramified cover of $\Gm^2$.
Theorem 1.2 admits the following consequence: \medskip

\noindent{\bf Theorem 1.3}. {\it Let $Y$ be a smooth affine surface with
log Kodaira dimension two. Suppose that there exists a cyclic (ramified) cover
$Y\rightarrow\Gm^2$. Let us embed  $Y$ in a projective space $\Pr_N$.
There exists a positive number $c$, depending on $Y$ and its projective embedding
such that  for every curve $\c\subset Y$,
$$
\deg(\c)\leq c\cdot \max\{1,\chi(\c)\}.
$$
}

\medskip

We observe that the condition that $Y$ is a ramified cover of $\Gm^2$
appears also in the main theorem of [NWY], where however the cyclicity condition
is not assumed and the torus $\Gm^2$ might be replaced by any arbitrary semi-abelian variety.
The authors prove in [NWY] that under the same hypothesis on the log Kodaira dimension every
entire curve on $Y$ is degenerate.
\medskip

We shall prove in \S 4 that no result like Theorem 1.1 can hold if the divisor
$D$ is a cubic instead of a quartic (in which case the
log canonical divisor of the complement would be zero). We shall also show in Prop. 4.4 that the hypothesis that the divisor $D$ has normal crossings cannot be dispensed with. \medskip

When the divisor $D$ is an {\it irreducible} cubic,  Beukers proved [Be, Thm. 3.3]
 in the arithmetic case,  that the integral points on $\Pr_2\setminus D$ are Zariski-dense, at least in a suitable extension of the ring of integers.  Beukers method also holds  over function fields, but in that case does not provide  automatically integral points of arbitrary high height.   \medskip

A possible generalization of Theorem 1.2 (so also of Theorems 1.1. and 1.3)
concerns the case of  polynomials $A(X,Y)\in\kappa(\cc)[X,Y]$, i.e. with non-constant
coefficients. Geometrically, this situation corresponds to the data of a three-fold $X$ and
a projection $\theta:X\rightarrow \c$. The solutions $(u_1,u_2,y)$ to (1.3) correspond
to sections of $\theta$. The method developed in the present paper permits to obtain the
analogue results in this setting too. 

 Another application of the present method concerns rational
curves on diagonal or generalized Fermat hypersurfaces.
For simplicity here  we limit ourselves  to the above statements, leaving such
generalizations to  a future paper.

\medskip

This paper is organized as follows: in the next section,  we shall prove a diophantine estimate (Theorem 2.1, Theorem 2.2 and its Corollary) concerning  the ``greatest common divisor" of two rational functions on $\cc$ of the form $a-1,b-1$, where $a,b$ have all their zeros and poles in a prescribed set.  This result, which is the function field analogue of Proposition 2 of [CZ2], is of independent interest and admits further applications which will be developed elsewhere.

In \S 3 we shall apply such estimates to  bound the degree of the solutions
of equations of the form $y^2=u_1^2+\lambda u_1+u_2+1$, where $u_1,u_2$ are, as above, rational
functions with prescribed zero and pole sets, and $\lambda$ is an arbitrary  scalar. As a
consequence, we shall obtain the inequality of Theorem 3.4 concerning the solutions to equation
(1.2).  Theorem 1.1, which is reformulated in equivalent form as Theorem 3.1,  is deduced from
Theorem 3.4.  The proof of Theorems 1.2 (which is a generalization of Theorem 3.4)
and 1.3 (generalizing Theorem 1.1)  runs along the same lines; a quick argument will be given at
the end of \S 3. In \S 4 we prove that the hypothesis of Theorem 1.1 cannot be weakened,  neither
by taking for the divisor $D$ a reducible curve of degree $\leq 3$  (instead of $4$), nor by
removing the hypothesis on the normal crossing singularities. Finally, in an appendix we shall
show some  arithmetic result in the spirit of Liardet-Raynaud theorems on torsion points on
curves.   \medskip

{\bf Sketch of the proofs}. For the reader's convenience, we now give here a brief outline
of the proof of Theorem 1.1. First of all an affine curve $\c$ in the surface
$X=\Pr_2\setminus D$ corresponds to a morphism from the normalization $\cc\setminus S$ of $\c$ to
$X$. Such a morphism, in turn, corresponds to the solution of the ``diophantine equation"  (1.2)
in rational functions $y,u_1,u_2\in\kappa(\cc)$ with $u_1,u_2$ $S$-units (see Lemma 3.3). We write
(1.2) as $y^2=A(X,Y)$ after setting  $A(X,Y)= X^2+\lambda X+Y+1$; differentiating both sides (with
respect to a suitable differential operator to be introduced  in Lemma 3.5) one gets another
equation satisfied by $y$, of the form  $2 y y^\prime=B(u_1,u_2)$  where the polynomial
$B(X,Y)=\uno X{\partial A\over \partial X}+\due Y{\partial A\over\partial Y}$ has its coefficients
in the field $\kappa(\cc)$. It turns out that the height of such coefficients is bounded in terms
of $S$ only, so in particular it is independent of the height of $u_1,u_2,y$.  Now we obtain that
the values of the two polynomials $A(X,Y), B(X,Y)$, calculated at the $S$-unit points $u_1,u_2$,
have ``many" zeros in common, i.e. the zeros of $y$. These facts can be reformulated in terms of
a functional greatest common divisor.

Let us consider for instance the particular case
$A(X,Y)=1+X+Y$ (so the equation (1.2) becomes $y^2=1+u_1+u_2$).  From
$2yy^\prime=u_1^\prime+u_2^\prime=\uno u_1+\due u_2$ one gets
$$
\uno y^2-2yy^\prime = \uno+\left(\uno-\due\right)u_2=
-\uno\left(u_2\left(\due-\uno\right){u_1\over u_1^\prime}-1\right)
=-\uno(w_2-1)
$$
where $w_2=u_2(\uno\due-1)$.
Also by symmetry
$$
\due y^2-2yy^\prime = -\due(w_1-1)
$$
with $w_1=u_1(\uno\due-1)$. The rational functions $w_1,w_2$ are $S^\prime$-units for a finite set
$S^\prime$ whose cardinality is bounded in terms of $\sharp(S)$  (and of course the chosen
derivation ${}^\prime$), but does not depend on $u_1,u_2$.
Then we are in the function-field analogue of the situation studied in [CZ2], where we proved a
bound for $\gcd(w_1-1,w_2-1)$ for units $w_1,w_2$ in a number field. (We note that this bound is
directly linked to Vojta's conjecture on diophantine approximation on certain surfaces, as
remarked by Silverman [S].)  As in the arithmetic case,  we obtain such a bound, actually in
stronger form; this is done in the next paragraph. Reinterpreting this inequality in terms of the
morphism corresponding to $(u_1,u_2,y)$, we obtain a bound for the degree of the image of such
morphisms, consisting in the inequality (1.1) of Theorem 1.1. \bigskip

\bigskip

\noindent{\bf \S 2. Greatest common divisors over function fields}.
 \smallskip

Let $\kappa$ be, as before, an algebraically closed field
of characteristic zero, $\cc$  a smooth complete curve
defined over $\kappa$, of genus $g=g(\cc)$,
$S\subset\cc$ a finite nonempty set of points of $\cc$.
We let $\O_S:=\kappa[\cc\setminus S]$ be the ring
of $S$-integers, i.e. of regular functions on the affine curve $\cc\setminus S$, and
$\O_S^*$ be the group of $S$-units, i.e. of rational functions
on $\cc$ with all their zeros and poles in $S$.
As usual, for a rational function $a\in\kappa(\cc)$,
we let $H(a)=H_{\cc}(a)$ be its height, i.e. its degree as a morphism
$a:\cc\rightarrow\Pr_1$. For $n\geq 2$ and
elements $u_1,\ldots,u_n\in\kappa(\cc)$ not all zero, we denote
by $H_{\cc}(u_1:\ldots:u_n)$, or simply $H(u_1:\ldots:u_n)$,  the projective height
$$
H(u_1:\ldots u_n)=
H_{\cc}(u_1:\ldots:u_n)=-\sum_{v\in \cc}\min\{v(u_1),\ldots,v(u_n)\}.
$$
Also, for $u\in \kappa(\cc)$, we let
$$
H_S(u)=-\sum_{v\not\in S}\min\{0,v(u)\}
$$
be the number of poles (with multiplicity) outside $S$.
Plainly $0\leq H_S(u)\leq H(u)$.
\medskip

As in the previous section, the Euler characteristic
$$
\chi=\chi(\cc\setminus S)=\sharp(S)+2g(\cc)-2
$$
will often appear.
We also note that if a nonconstant $S$-unit exists, as in our
context, then necessarily $\sharp(S)\geq 2$, so $\chi\geq 0$.

We also make the following simple but important remark: if we replace the curve by a cover of it,
of degree $d$,  then the new heights will be multiplied by $d$ while the new $\chi$ will be at
least the old one multiplied by $d$. (This follows at once from the Riemann-Hurwitz formula, taking
into account the possible ramifications in  and out of $S$.) By this observation it will be clear
that in proving the statements in this section we may often consider the minimal function field
containing the relevant quantities. \medskip

\noindent{\bf Proposition 2.1.}\ {\it Let $a,b\in \O_S^*$ be
multiplicatively independent $S$-units, not both constant, and let
$h,k$ be positive integers. Then, either $$H(a)\le
h\cdot[\kappa(\cc):\kappa(a,b)],\qquad  H(b)\le
k\cdot[\kappa(\cc):\kappa(a,b)]$$ or $$ H_S\left({1-a\over
1-b}\right)\ge {hk\over hk+h+k}H(b)-{k\over
hk+h+k}(H(a)+H(b))-{hk+h+k-1\over 2}\chi. $$}\medskip

Another, essentially equivalent,  formulation of the result involves a kind of
``$S-\gcd$" of
$1-a,1-b$, which here we measure by counting the number of common zeros outside $S$.
Namely:\medskip

\noindent{\bf  Proposition 2.2.}\ {\it Let $a,b\in \O_S^*$ be
multiplicatively independent $S$-units, not both constant, and let
$h,k$ be positive integers. Then, either $$H(a)\le
h\cdot[\kappa(\cc):\kappa(a,b)],\qquad  H(b)\le
k\cdot[\kappa(\cc):\kappa(a,b)]$$ or $$ \sum_{v\not\in
S}\min\{v(1-a),v(1-b)\}\le {h+2k\over hk+h+k}H(b)+{k\over
hk+h+k}H(a)+{hk+h+k-1\over 2}\chi. $$}\medskip

For the sake of generality, we  have let here $h,k$ be any positive integers,
but special choices of them obviously may lead to simpler bounds;
e.g., given $\epsilon>0$, on choosing $h=k\sim 4\epsilon^{-1}$
 we see  that the right side of  the displayed inequality in Proposition
$2.2$ becomes $\le\epsilon(H(a)+H(b))+O_{\sharp S,g,\epsilon}(1)$.
More explicitly, we have for instance the following corollary,
in which we also add a complementary part concerning the case of
multiplicatively dependent $a,b$.
In this  case we have some nontrivial relation
$a^r=\lambda b^s$ for integers $r,s$ not both zero and a
$\lambda\in \kappa^*$; if not both $a,b$ are constant,
such pairs $(r,s)$ form a $1$-dimensional  lattice in $\Z^2$ and if $(r,s)$
 generates the lattice we speak of a generating relation. We have:\medskip

\noindent{\bf Corollary 2.3.}\ {\it Let $a,b\in \O_S^*$ be $S$-units, not
both constant, and let   $H:=\max\{H(a),H(b)\}$.

(i) If $a,b$ are multiplicatively independent, we have $$
\sum_{v\not\in S}\min\{v(1-a),v(1-b)\}\le  3{\root 3\of
2}(H(a)H(b)\chi)^{1\over 3}\le 3{\root 3\of 2}(H^2\chi)^{1\over
3}. $$

(ii) If $a,b$ are multiplicatively dependent, let $a^r=\lambda b^s$ be a
generating relation. Then either $\lambda\neq 1$ and
$\sum_{v\not\in S}\min\{v(1-a),v(1-b)\}=0$, or $\lambda=1$ and
$$
\sum_{v\not\in S}\min\{v(1-a),v(1-b)\}\le \min\left\{
{H(a)\over|s|},{H(b)\over |r|}\right\}\le {H\over\max\{|r|,|s|\}}.
$$
}\medskip

The constant $3{\root 3\of 2}$ in part (i) cannot be replaced by
anything smaller than $(4/3)^{1/3}$, as shown by the example
$a=t^3$, $b=-t(t+1)$.

It is amusing that this corollary, which concerns function
fields and is proved here
 by function field-methods, in turn easily implies the (well-known) arithmetical
conclusion of the finiteness of torsion points on an irreducible  curve in
$\Gm^n$ which is not a translate of an algebraic subgroup.
In an Appendix we shall briefly show this deduction, and also some relations of the results
with papers by Ailon-Rudnick [AR] (on $\gcd(f^n-1,g^n-1)$ for polynomials $f,g$)
and Bombieri-Masser-Zannier [BMZ] (on the intersections of a given curve in $\Gm^n$
with the algebraic subgroups).\medskip

\noindent{\bf Proof of Proposition 2.1}. We begin with some  preliminary observations.

We have already noted that we may increase the function field to a
finite extension without loss of generality. Therefore (since
$a,b$ are not both constant) we may assume that
$\kappa(\cc)=\kappa(a,b)$, namely we may assume that $a,b$
generate the function field.

 Since $a,b$ are not both constant, $\sharp S\ge 2$, so $\chi\ge 0$;
hence if $b\in \kappa$  the result holds trivially. Let us then suppose that  $b\not\in \kappa$.

Let now $c=1-b$, so $b+c=1$, and let $S'$ be the set of zeros of $c$ outside $S$.
Since $b,c$ are nonconstant $(S\cup S')$-units, by Mason's
$abc$-theorem for function fields
(e.g., Cor. 1 of Thm. A of [BrM] with $n=3$) we
have $H(c)\le H(1:b:c)\le \sharp S+\sharp S'+2g-2$. But $H(c)$ is the number of zeros of $c$,
counted with multiplicity, so in particular $\sum_{v(c)>0}v(c)\le
\sharp S+\sharp S'+2g-2$. Since $v(c)>0$ for all $v\in S'$, we have
$$
\sum_{v(1-b)>0}(v(1-b)-1)\le \chi.\eqno(2.1)
$$\medskip

Next,  suppose we have a nontrivial relation $a^rb^s=\lambda\in \kappa^*$.
If $a-1$ and $b-1$  have no common zero, we have
$$
H_S\left({1-a\over 1-b}\right)\ge
\sum_{v(1-b)>0, v\not\in S}v(1-b)=H(1-b)-\sum_{v(1-b)>0,v\in S}v(1-b).
$$
The first term on the right equals $H(b)$; also,
$0\le \sum_{v(1-b)>0,v\in S}v(1-b)\le \sharp S-2+\sum_{v(1-b)>1,v\in S}(v(1-b)-1)$. Hence
by (2.1) we have $H_S({1-a\over 1-b})\ge H(b)-(2\sharp S+2g-4)$, proving the proposition and more.
Hence we may assume that $1-a,1-b$  have some common zero  $P$.
Then $a(P)=b(P)=1$, which yields $\lambda=1$,
against the multiplicative independence assumption.

Therefore in what follows we shall assume that $a,b$ are multiplicatively independent modulo
$\kappa^*$.\medskip

The argument now  will mimic [CZ2], but will involve a simple  Wronskian argument  rather than
the Schmidt Subspace Theorem.  (Especially after [BrMa], Wronskians are a
familiar tool in diophantine questions over  function fields.
See also [W] for an effective proof  with Wronskians of a function field version
 of the Subspace Theorem.) Thus we   recall some standard facts about Wronskians.

For elements $f_1,\ldots ,f_n\in \kappa(\cc)$ and a nonconstant $t\in \kappa(\cc)\setminus\kappa$ we let
the ``Wronskian" $W_t(f_1,\ldots ,f_n)$ be the  determinant of the $n\times n$
matrix whose $j$-th row-entries are the $(j-1)$-th derivatives of the $f_i$'s
with respect to $t$.
It is well known that $W_t=0$ if and only if the $f_i$'s are linearly dependent over
$\kappa$. (Recall that here char $\kappa=0$.)
Let $z\in K$ be another nonconstant element. Then we have the known, easily proved, formula
$$
W_z(f_1,\ldots ,f_n)=({dt\over dz})^{n\choose 2}W_t(f_1,\ldots ,f_n).\eqno(2.2)
$$
For a place $v\in\cc$ we choose once and for all a local parameter $t_v$ at
$v$ and we define $W_v:=W_{t_v}$. This depends on the
choice of $t_v$, but (2.2) shows that the order $v(W_v)$ depends only on $v$.

To prove the Proposition  we shall consider suitable Wronskians. \medskip

We  define $q=(1-a)/(1-b)$ and,  letting $n:=hk+h+k$, we define functions
$f_1,\ldots ,f_n$ as follows. For $i=1,\ldots ,k$ we let $f_i:=a^{i-1}q$
while we define $f_{k+1},\ldots ,f_n$ as the functions $a^rb^s$,
$r=0,1,\ldots ,k$, $s=0,1,\ldots ,h-1$, in some order.

We now choose   nonconstant $t$, $t_v$ as above  and we put
$$
\omega=W_t(f_1,\ldots ,f_n),\qquad  \omega_v=W_{t_v}(f_1,\ldots ,f_n).
$$

Suppose  first that $\omega=0$; then, as we have remarked, the $f_i$ are linearly dependent
over the constant field $\kappa$; recalling the definition of the $f_i$ this amounts to
a relation $P_1(a)(1-a)+P_2(a,b)(1-b)=0$, where $P_1(X)$ is a polynomial of degree
$\le k-1$ and $P_2(X,Y)$ is a polynomial of degree $\le k$ in $X$ and $\le h-1$
in $Y$ and where not both $P_1,P_2$ are  zero. Observe that $P_1(X)(1-X)+P_2(X,Y)(1-Y)$
is not identically zero, for otherwise $P_1(X)$ would vanish (set $Y=1$) and then
$P_2$ would also vanish, a contradiction.
Then  we  find a nontrivial polynomial relation $P(a,b)=0$,
for a polynomial $P\neq 0$ of degree $\le k$ in $X$ and $\le h$ in $Y$. Of course we may assume
that $P$ is irreducible.

Then, since $a,b$ generate our function field, we have $H(a)=\deg a\le h$, $H(b)=\deg b\le k$,
falling into the first possibility of the sought conclusion.

Therefore in what follows we assume that $\omega\neq 0$.\medskip

To go on, we first seek for a lower bound for $v(\omega_v)$ and for this we shall
distinguish several cases.  \medskip

Case (i): $v\not\in S$, $v(q)<0$. We first make a simple observation: suppose that
two functions $f_i,f_j$, $i\neq j$, have a pole at $v$ of the same order $>0$. Then
by subtracting from the $i$-th column a constant multiple of the $j$-th column, we
may suppose in calculating the Wronskian that
$f_i$ has a pole of smaller order than $f_j$ at $v$. Therefore, by repeating this
procedure we may assume that the functions $f_i$ which have a pole  at $v$ have in
fact poles of pairwise distinct orders, not exceeding the original maximal order.

In the present Case (i) the only functions which (may) have a pole at $v$ are
$f_1,\ldots ,f_k$, because $v\not\in S$, so $a$ and the remaining $f_i$'s are
units at $v$ by assumption.
 Hence,  the only  poles we shall possibly find in
$\omega_v$ will come from the first $k$
columns. However by the
above observation we may change the actual $f_i$'s, $i=1,\ldots ,k$, to assume that
the negative ones among
$v(f_1),\ldots ,v(f_k)$ are all distinct and $\ge v(q)$. Suppose that after such
column operations and suitable renumbering only $f_1,\ldots, f_r$ have a pole at $v$
and
$v(f_1)<\ldots <v(f_r)<0$. Plainly we will have $v(q)\le v(f_1)$ and $r\le k$, so
$r\le\min(k,-v(q))$. Also, observe that each derivation with respect to
$t_v$ increases the order of a pole by
$1$ and leaves regular a regular function at
$v$. In conclusion, by looking at the individual terms obtained in the expansion of
the determinant (after having performed the column operations), a simple inspection
shows that   we have $v(\omega_v)\ge v(f_1)+\ldots +v(f_r)+{r\choose 2}-r(n-1)$,
whence, since we are assuming the $v(f_i)$ to be distinct,
$$
v(\omega_v)\ge
rv(q)-r(n-r)\ge rv(q)+v(q)(n-r)=v(q)n.
$$
\medskip

Case (ii): $v\not\in S$, $v(q)\ge 0$. Now every element of the local
Wronskian matrix is
$v$-integral, so the same holds for the determinant, i.e., $v(\omega_v)\ge 0$ in this case.
\medskip

Case (iii): $v\in S$, $v(b)>0$. This case contains the crucial point.
As in [CZ2], we consider the identity
$$
a^jq-a^j(1-a)(1+b+\ldots +b^{h-1})=a^jb^hq.\eqno(2.3)
$$
This will be useful to approximate $a^jq$ with a polynomial in $a,b$, at the places
under consideration.

In fact, we may use the identity  to replace, for
$i=1,\ldots ,k$, the function $f_i$ with the left side of (2.3), with $j=i-1$,
which by (2.3) equals $a^{i-1}b^hq$, denoted  $g_i$. Observe that this
corresponds to subtract from
$f_i$ a certain $\kappa$-linear combination of $f_{k+1},\ldots ,f_n$,
and thus the value of $\omega_v$ is unchanged. We have $v(g_i)=(i-1)v(a)+hv(b)+v(q)$.
Since $v({d^lf\over dt_v^l})\ge v(f)-l$,
we easily find, on looking again at the individual terms in
the determinant expansion, that
$$
v(\omega_v)\ge {k(k-1)\over
2}v(a)+hkv(b)+kv(q)+\left(\sum_{i=k+1}^nv(f_i)\right)-{n\choose 2}.
$$\medskip

Case (iv): $v\in S$, $v(b)\le 0$. We now argue directly with the terms in the
determinant expansion (that is, we do not perform any column operation). Since
$v(f_i)=(i-1)v(a)+v(q)$ for $i=1,\ldots ,k$, we find as in the previous case that
$$
v(\omega_v)\ge {k(k-1)\over 2}v(a)+kv(q)
+\left(\sum_{i=k+1}^nv(f_i)\right)-{n\choose 2}.
$$\medskip

Summing over all places $v$ of $\kappa(\cc)$, taking into account the estimates obtained in the four cases, and recalling that $\sum_{v\in S}v(a)=\sum_{v\in S}v(b)=\sum_{v\in S}v(f_i)=0$ for $i>k$, because $a,b$ are $S$-units, we thus get
$$
\sum_vv(\omega_v)\ge \sum_{v\not\in S, v(q)<0}nv(q)
+hk\sum_{v\in S,v(b)>0}v(b)+k\sum_{v\in S}v(q)-{n\choose 2}\sharp S.\eqno(2.4)
$$

Now, $\sum_{v\in S,v(b)>0}v(b)=\sum_{v(b)>0}v(b)$, because $b$  is an $S$-unit; also,
$\sum_{v(b)>0}v(b)=-\sum_{v(b)<0}v(b)=H(b)$. Moreover,
(2.2) with $z=t_v$ shows that $v(\omega_v)={n\choose 2}v(dt/dt_v)+v(\omega)$.
On summing over $v$ this yields
$$
\sum_vv(\omega_v)={n\choose 2}\sum_vv({dt\over dt_v})+\sum_vv(\omega)=
{n\choose 2}(2g-2),
$$
the last equality holding because of the product formula
(for $\omega\in \kappa(\cc)^*$) and because $2g-2$ is the degree of any canonical divisor.
Comparing with the above yields
$$
{n\choose 2}\chi- \sum_{v\not\in S, v(q)<0}nv(q)
\ge hkH(b)+k\sum_{v\in S}v(q).\eqno(2.5)
$$
Finally, $-\sum_{v\in S}v(q)\le H(q)\le  H(a)+H(b)$, whence
$$
H_S(q)\ge {hk\over hk+h+k}H(b)-{k\over hk+h+k}(H(a)+H(b))
-{n-1\over 2}\chi,
$$
concluding the proof.\medskip

\noindent{\bf Proof of Proposition $2.2$.} Note that by definition
$H_S({1-a\over 1-b})=\sum_{v\not\in S, v(1-b)>v(1-a)}(v(1-b)-v(1-a))$.
In turn the right side may be plainly replaced by $\sum_{v\not\in S}(v(1-b)-\min(v(1-b),v(1-a)))$.

Now, $\sum_{v\not\in S}v(1-b)$  does not exceed $H(1-b)=H(b)$, whence
$$
H_S({1-a\over 1-b})\le H(b)-\sum_{v\not\in S}\min\{v(1-a),v(1-b)\}.
$$
An application of Proposition 2.1 now immediately leads to the sought
inequality.\medskip

\noindent{\bf Proof of Corollary $2.3$.}  If one of $a,b$ is constant, then it cannot be $1$
because $a,b$  are multiplicatively independent. Therefore the sum vanishes and we are done. So,
let us suppose that none of $a,b$ is constant. Then we may also assume that $H(a)\ge
H(b)>0$.\medskip

Let us first deal with part (i),  supposing  that
$\chi=0$ to start with; then $\#S=2,g=0$ and necessarily there is a relation
$a^rb^s=\gamma\in\kappa^*$, where $r,s$
 are not both zero: this is because some function $a^{\deg b}b^{\pm\deg a}$ has no zeros or poles.
 Then $\gamma\neq 1$ by assumption, so $\min(v(1-a),v(1-b))\le 0$ for all $v$, proving the
result.

Therefore, suppose
$\chi\neq 0$ in  the sequel.

Also, as we have remarked, we may assume that $\kappa(\cc)=\kappa(a,b)$.\medskip

Let us choose now $h=[(4H(a)^2/H(b)\chi)^{1/3}]-1$,  $k=[(4H(b)^2/H(a)\chi)^{1/3}]-1$.

Suppose that $k<1$. Then $[(4H(b)^2/H(a)\chi)^{1/3}]<2$, whence $(4H(b)^2/H(a)\chi)^{1/3}<2$ and
so $H^2(b)<2H(a)\chi$. Then $H^3(b)<2H(a)H(b)\chi$.

In this case we use the obvious inequality
$\sum_{v\not\in S}\min(v(1-a),v(1-b))\le H(1-b)=H(b)<(2H(a)H(b)\chi)^{1/3}$, proving what we
need, with a better constant.\medskip

Hence in what follows we assume that  $h\ge k\ge 1$, so in particular we may apply Proposition
2.2. The conclusion gives two possibilities.

In the first case we have $H(a)\le h$, whence $H(a)^3\le (h+1)^3\le 4H(a)^2/H(b)\chi$, whence
$H(a)H(b)\chi\le 4$. Since $\chi\ge 1$, this implies $H(b)\le 2$ so as above we obtain
$\sum_{v\not\in S}\min(v(1-a),v(1-b))\le H(1-b)=H(b)\le 2$. Again this gives the sought result
since $H(a)H(b)\chi\ge 1$.\medskip

Then we may assume $H(a)>h$, so the second alternative of Proposition 2.2 must hold.

It is easily checked that the coefficient $(h+2k)/(hk+h+k)$ is bounded by $3/(k+2)$, since $h\ge
k$. Similarly, $k/(hk+h+k)\le 1/(h+2)$. Therefore
$$
\sum_{v\not\in S}\min(v(1-a),v(1-b))\le {3\over k+2}H(b)+{1\over h+2}H(a)+{(h+1)(k+1)-2\over
2}\chi.
 $$

We now use $k+2\ge (4H(b)^2/H(a)\chi)^{1/3}$, $h+2\ge
(4H(a)^2/H(b)\chi)^{1/3}$ and also $(h+1)(k+1)\le
(4H(b)^2/H(a)\chi)^{1/3}(4H(a)^2/H(b)\chi)^{1/3}=(16H(a)H(b)/\chi^2)^{1/3}$,
which yields $$ \sum_{v\not\in S}\min(v(1-a),v(1-b))\le
4({H(a)H(b)\chi\over 4})^{1/3}+ {1\over
2}(16H(a)H(b)\chi)^{1/3}=3{\root 3\of 2}(H(a)H(b)\chi)^{1/3}
$$ concluding the proof of part (i).\medskip

For part (ii), if $\lambda\neq 1$ the result follows as in part (i), so
suppose $a^r=b^s$ in a generating relation. Plainly we have that
$\gcd(r,s)=1$ whence $hr+ks=1$ for some integers $h,k$. Then, setting
$c:=a^kb^h$, which is an $S$-unit,  we have $a=c^s,b=c^r$. This yields $\min
(v(1-a),v(1-b))=v(1-c)+\min(v({1-c^s\over 1-c}),v({1-c^r\over 1-c}))$. This
last term vanishes for $v\not \in S$, since then $c$ is regular at $v$ and the
polynomials ${1-X^s\over 1-X}$ and ${1-X^r\over 1-X}$ are coprime.
Therefore the relevant sum equals $\sum_{v\not\in S}v(1-c)\le H(c)$ and the result follows
at once since $|s|H(c)=H(a)$, $|r|H(c)=H(b)$. \medskip

\bigskip


\noindent{\bf \S 3.  Maps to the complement of a conic and two lines}.
\smallskip

Let us consider the configuration of a conic and two lines in $\Pr_2$:
let $D_1$ be a smooth conic in $\Pr_2$, $D_2,D_3$ be
distinct lines, all defined over $\kappa$. We are particularly interested
in the case where $D_1,D_2,D_3$ are in general position, by this
meaning that the two lines $D_2,D_3$ intersect  the conic $D_1$
at four distinct points. Such configurations  form a one-dimensional family,
parametrized by the cross-ratio of the four intersection points
in $D_1\cap(D_2\cup D_3)$ with respect to the conic.

Let now $\cc$ be a smooth complete algebraic curve of genus $g(\cc)$
defined over $\kappa$,  and  $S\subset \cc$ be a finite nonempty set of points of $\cc$.
 We shall denote again by $\O_S$ the ring $\kappa[\cc\setminus S]$ of regular functions on the affine
curve $\cc\setminus S$; its elements will also be called $S$-integers; then $\O_S^*$ will be its
group of units,  i.e. the multiplicative group of rational functions on $\cc$ having all their
zeros and poles in $S$;  they will be called $S$-units. 
We are interested in classifying regular maps from $\cc\setminus S$
to $\Pr_2$  ``omitting" the divisor $D:=D_1+D_2+D_3$,
i.e. morphisms $f:\cc\rightarrow \Pr_2$ such that $f^{-1}(D)\subset S$.

\smallskip

Before passing to the ``general position" situation, i.e. when the lines
$D_2,D_3$ intersect the conic $D_1$ at four distinct points we analyse
some ``degenerate" cases  which   deserve attention.
\smallskip

A first interesting case comes from a tangent line $D_2$ to the conic $D_1$ and a line $D_3$ with
$D_1\cap D_2\cap D_3=\emptyset$. This case is easily reduced to the ``diophantine equation"
$$
1+u_1+u_2=y^2\eqno(3.1)
$$
to be solved in units $u_1,u_2\in\O_S^*$ and regular functions $y\in \O_S$.

This situation can also be recovered  from the ``general position"  case
where the four intersection points  of $(D_2\cup D_3)\cap D_1$
 have cross-ratio $-1$ with respect to the conic $D_1$ (the reduction can
be obtained via the unramified covering $v^2=u_2$).
The special case where moreover $\cc=\Pr_1$ has genus zero and $S$ has three points
has been treated by the second author in [Z2], using different tools.
\smallskip

A second degenerate case is provided by two distinct non tangent lines
$D_2,D_3$ to the conic $D_1$  such that $D_1\cap D_2\cap D_3$ is non empty
(so necessarely consists of a single point).
The problem in this case can be reduced to the equation
$$
y={u_1-1\over u_2-1}
$$
to be solved in $S$-units $u_1,u_2$ and $S$-integers
$y\in\O_S$. So it consists of a pure divisibility problem in the
ring $\O_S$. In this case, however, no bound of the form (1.1) holds
on the corresponding surface $\Pr_2\setminus D$ (this is the content of Proposition 4.3).
 \medskip

From now on we shall suppose that $D_1,D_2,D_3$ are in general position,
and shall denote by $D$ also the support of the divisor $D_1+D_2+D_3$.
We begin by restating Theorem 1.1 in the following way.
\medskip

\noindent{\bf Theorem 3.1}. {\it Let $\cc,S,D$ be as above. Let
$f:\cc\rightarrow\Pr_2$ be a
non constant morphism such that $f^{-1}(D)\subset S$. Then the degree of
the curve $f(\cc)$ verifies
$$
\deg(f(\cc))\leq 2^{15}\cdot 35\cdot \max\{1,\chi(\cc\setminus S)\}. \eqno(3.2)
$$
}

\bigskip

Let us introduce the {\bf notation for the proofs}.
As in the previous paragraph,
we associate with each point $v$ of $\cc$  a discrete valuation
of the function field $\kappa(\cc)$, trivial on $\kappa$,
normalized so that its value  group is the group of integers $\Z$
and shall denote it  by the same letter $v$.
The height of a rational function  $a\in \kappa(\cc)$
will be denoted now
by the symbol $H_{\cc}(a)$, since the reference to the curve
$\cc$ will be relevant in the sequel; recall that it
coincides with the degree of $a$ viewed as a morphism
$\cc\rightarrow\Pr_1$;  hence it is given by the formula
$$
H_\cc(a)=\sum_{v\in\cc}\max\{0,v(a)\}.
$$
If $\dd\rightarrow\cc$ is a dominant morphism of smooth (irreducible) projective curves,
then the function field of $\cc$ injects into the function field of $\dd$;
in that case an element $a\in\kappa(\cc)$ could also be viewed as a rational function on $\dd$.
We shall then write $H_{\dd}(a)$ to denote its degree as a function on $\dd$; it verifies 
$H_{\dd}(a)=[\kappa(\dd):\kappa(\cc)]\cdot H_{\cc}(a)$.
For every element $a\in \kappa(\cc)$, we write $d(a)$ for its differential.
Finally, recall that the Euler characteristic
of the affine curve $\cc\setminus S$ is by definition the integer
$ \chi_S(\cc)=\chi(\cc\setminus S)=2g(\cc)-2+\sharp(S).$
\smallskip

We begin by proving the following
\medskip

\noindent {\bf Lemma 3.2}. {\it Let $\cc,S$ be as before, and
$f:\cc\setminus S\rightarrow\A^2$
be a regular map from the affine curve $\cc\setminus S$ to the affine plane,
given as $f(p)=(x(p),y(p))$ for $S$-integers $x,y\in \O_S$.
Then the degree of the image $f(\cc)$
(viewed as embedded canonically in $\Pr_2$) verifies}
$$
\deg(f(\cc))\leq H_{\cc}(x)+H_{\cc}(y).
$$
\smallskip

\noindent {\it Proof}. There exists a choice of scalars $(a,b,c)\in \kappa^3$ with
$(a,b)\neq (0,0)$ such that the line of equation $ax+by=c$ does not intersect the curve
$f(\c)$ (where as usual $\c:=\cc\setminus S$) at infinity nor at any singular point  and is not
tangent at any point (actually this is true for a generic choice of $a,b,c$). Then the degree of
$f(\cc)$ is  the number of intersection points of the affine curve $f(\c)$  with such a
line. This  number is bounded by the number of solutions in $\cc\setminus S$ to the equation
$ax(p)+by(p)=c$; now the (regular) function $ax(p)+by(p)$ has height
$\leq H_\cc(x)+H_\cc(y)$, so the estimate follows.\cvd

\medskip

The next lemma will rely the existence of morphisms
$\cc\setminus S \rightarrow X$ to the solutions of a certain diophantine
 equation over the function field $\kappa(\cc)$.
\medskip

\noindent {\bf Lemma 3.3}. {\it Suppose
$f:\cc\setminus S\rightarrow \Pr_2\setminus D$  is a morphism.
There exist a scalar $\lambda\in \kappa$, $S$-units $u_1,u_2\in\O_S^*$
and an $S$-integer $y\in\O_S$ satisfying
$$
H_{\cc}(y)+H_{\cc}(u_1)\geq \deg(f(\cc))
$$
such that
$$
y^2=u_1^2+\lambda u_1+u_2+1.\eqno(3.3)
$$
}
\medskip

\noindent{\it Proof}.
We can certainly choose homogeneous coordinates $(x_0:x_1:x_2)$
in $\Pr_2$ such that the intersection
$D_1\cap D_2$ consists of the two points $(0:1:1),(0:1:-1)$, while the intersection
$D_1\cap D_3$ is formed by the points $(1:0:1), (1:0:-1)$.
Then the lines $D_2,D_3$ are defined  by the equations
$$
D_2:\ x_0=0\qquad {\rm and}\qquad D_3:\ x_1=0
$$
respectively. Also the conic $D_1$ must belong to the pencil of conics
 defined by the equation
$$
x_2^2-x_1^2-\lambda x_0x_1-x_0^2=0,\eqno(3.4)
$$
for a suitable $\lambda\in\kappa$.
We shall now put $x:=x_1/x_0,\, y=x_2/x_0$, which are affine coordinates
for the plane $\Pr_2\setminus D_2\simeq {\bf A}^2.$ Note that the line $D_3$
has the equation $x=0$ in such affine coordinates.
Let now $f:\cc\setminus S \rightarrow\Pr_2\setminus D$ be a regular map.
The fact that $f(\cc\setminus S)$ avoids the line at infinity $D_2$ and the line $D_3$
means precisely that $f$ can be written in affine coordinates as
$$
f(p)=(u_1(p),y(p))\qquad (p\in \c\setminus S)
$$
for a regular function $y\in\O_S$ and a unit $u_1\in\O_S^*$.
Imposing the condition that the image of $f$ avoids also the conic
defined by (3.4)  means that the regular function
$$
u_2:=y^2-u_1^2-\lambda u_1-1
$$
is in fact a unit. Finally, the bound $H_{\cc}(y)+H_{\cc}(u_1)\geq \deg(f(\cc))$
follows from the preceding Lemma, concluding the proof.
\cvd

\medskip

We shall now work with equation (3.3). Our aim is to prove the following
\medskip

\noindent{\bf Theorem 3.4}. {\it Let $\cc,S$ and $\lambda$ be as before;
then every solution $(y,u_1,u_2)\in\O_S\times(\O_S^*)^2$
of equation (3.3) satisfies one of the following conditions:
\smallskip

\item{(i)} a sub-sum on the right term of (3.3) vanishes;
\smallskip

\item{(ii)} $u_1,u_2$ verify a multiplicative dependence relation of the form
$u_1^r\cdot u_2^s=\mu$, where $\mu\in\kappa^*$ is a scalar and $r,s$ are
integers, not both zeros, and both $\leq 5$;

\smallskip
\item{(iii)} $\max \{H_{\cc}(u_1),H_{\cc}(u_2)\}\leq 2^{14}\cdot 35\cdot \chi_S(\cc)$.

}
\medskip

The idea of the proof is that if the right-side term in (3.3) is a square, then it has
``many" zeros in common with  its ``derivative", which can also be written
as a linear combination (with non-constant coefficients) of $S$-units.
To exploit this remark we should first specify the precise meaning of ``derivative";
for this reason we state the following \medskip

\noindent{\bf Lemma 3.5}. {\it There exists a differential form $\omega$ on
$\cc$ and a finite set $T\subset\cc$ of cardinality $\sharp(T)=
\max\{0,2g(\cc) -2\}$ such that
for every  $S$-unit $u\in\O_S^*$ there exists an $(S\cup T)$-integer
$\theta_u\in \O_{S\cup T}$ having only simple poles such that
$$
{d(u)\over u}=\theta_u\cdot\omega,\qquad H_{\cc}(\theta_u)\leq
 \chi_S(\cc).\eqno(3.5)
$$
}

\smallskip

\noindent{\it Proof}. Suppose first the genus $g(\cc)\neq 0$.
Then there exists a regular differential form $\omega$,
having $2g(\cc)-2$ distinct zeros of multiplicity one. Let $T$ be
its zero set, which has naturally cardinality $ 2g(\cc)-2$.
Let now $u\in\O_S^*$ be an $S$-unit.
Then the only rational function $\theta_u$ on $\cc$ satisfying (3.5)
can have poles only at the zeros of $\omega$ or the poles of $d(u)/u$,
hence it is an $(S\cup T)$ integer. The poles of $d(u)/u$ are
necessarily contained in $S$ and of multiplicity $\leq 1$.
It follows that the number of poles with multiplicity of $\theta_u$ is
$\leq\sharp(S)+2g(\cc)-2$ as wanted.
Consider now the case of genus $0$.
Let $v,w$ be distinct points on $S$ (note that $\sharp(S)\geq 2$,
otherwise all $S$-units would be constant) and let $\omega$
be a differential form having two simple poles in $v,w$ and no zero.
Then the rational function $\theta_u$ appearing in (3.5) can have poles
only in $S\setminus \{v,w\}$, of multiplicity at most one, so its height is again
$\leq \sharp(S)-2=\chi_S(\cc)$ as wanted.\cvd
\medskip

We shall now fix such a differential form $\omega$ as in Lemma 3.5,
and also the finite set of places $T$ (which will be empty if the genus
of $\cc$ is zero or one,  and  of cardinality $2g(\cc)-2$
for $g(\cc)\geq 2$). For a rational function $a\in \kappa(\cc)$,
we shall denote by $a^\prime$ the only rational function such that
$$
d(a)=a^\prime\cdot \omega;
$$
in such a notation, the rational function $\theta_u$ appearing in the
previous lemma equals $u^\prime/u$.
The analogous estimate for $S$-integers (instead of $S$-units) is
\medskip

\noindent{\bf Lemma 3.6}. {\it Let $a\in\O_S$ be an $S$-integer.
Then we have
$$
d(a)=a^\prime\cdot\omega
$$
where $a^\prime$ is a $(S\cup T)$-integer whose height  satisfies
$$
H_\cc(a^\prime)\leq H_\cc(a)+\chi_S(\cc).
$$
Moreover for every place of $v\in S\cup T$, we have $v(a^\prime)\geq v(a)-1$.
}

\medskip

\noindent{\it Proof}.
The proof is analogous to the previous one. As before,
we first treat the case of positive genus. The poles of $a^\prime$ are
either the zeros of $\omega$ (and we have $2g(\cc)-2$ of them counting multiplicities)
or the poles of $a$.
If $v$ is such a pole, then $v\in S$ since $a$ is supposed to be
an $S$-integer, and  $v(a^\prime)=v(a)-1$ (if $v$ is not
a zero of $\omega$). Then the difference between the number
of poles with multiplicities of $a^\prime$ and $a$ is bounded by
$\sharp(S)+2g(\cc)-2$ as wanted.
The genus zero case  is treated in the same way as before.\cvd

\medskip

Our next goal will be the proof of a formal identity for the derivative of the values of a
polynomial \smallskip

\noindent{\bf Lemma 3.7}. {\it Let $A(X,Y)\in\kappa[X,Y]$ be a polynomial,
$u_1,u_2$ be $S$-units in $\kappa(\cc)$; let $B(X,Y)\in\O_{S\cup T}[X,Y]$
be the polynomial defined by
$$
B(X,Y)=\uno\cdot X{\partial A\over \partial X}(X,Y)+
\due\cdot Y{\partial A\over\partial Y}(X,Y).
$$
Then, with the above notation for the $1$-form $\omega$,
the differential of $A(u_1,u_2)$ is
$$
d(A(u_1,u_2))=B(u_1,u_2)\cdot\omega.
$$
}
\smallskip

Note that by using the symbol ${}^\prime$ as before we can rewrite the above identity as
$(A(u_1,u_2))^\prime=B(u_1,u_2)$.
\smallskip

\noindent{\it Proof}. Write $A(X,Y)=\sum_{i,j} a_{i,j}X^iY^j$, where
the sum runs over a finite subset of $\N^2$, and the coefficients
$a_{i,j}$ are scalars in $\kappa$. The differential of each monomial
calculated in $(u_1,u_2)$ is
$$
d(a_{i,j}u_1^i u_2^j)=\left(a_{i,j}\cdot i\cdot u_1^\prime u_1^{i-1}u_2^j+
a_{i,j}\cdot j\cdot u_1^{i}u_2^\prime u_2^{j-1}\right)\cdot \omega=
a_{i,j}\left(\uno\cdot i\cdot u_1^iu_2^j+
\due\cdot j\cdot  u_2^\prime u_1^{i}u_2^j\right) \cdot\omega
$$
and the lemma follows.
\cvd

We now return to equation (3.3).
\medskip

\noindent{\bf Lemma 3.8}. {\it Let $A(X,Y)\in\kappa[X,Y]$ and
$B(X,Y)\in\O_{S\cup T}(\cc)[X,Y]$ be  
$$
\eqalign{
A(X,Y)&=X^2+\lambda X+ Y+1\cr
B(X,Y)&=2\uno\cdot X^2+\lambda\uno\cdot X+\due\cdot Y;
}\eqno(3.6)
$$
let $F(X)\in\O_{S\cup T}[X], G(Y)\in\O_{S\cup T}[Y]$ be the resultants of $A(X,Y), B(X,Y)$ with
respect to $Y$ and $X$, i.e. the polynomials
$$
\eqalign{
F(X)&:=\left(\due-\uno\right)X^2+\lambda\left(\due-1\right)X+\due=
{\rm Res}_Y(A(X,Y),B(X,Y)),\cr
G(Y)&:=2\left(2\uno-\due\right)Y^2+
\left((\lambda^2-2)\uno\due+(4-\lambda^2)\uno+4\left(\uno\right)^2\right)Y+
2(4-\lambda^2)\left(\uno\right)^2\cr
&={\rm Res}_X(A(X,Y),B(X,Y)).
}\eqno(3.7)
$$
For every solution $(y,u_1,u_2)\in\O_S\times(\O_S^*)^2$ of (3.3)
we have
$$
\eqalign{
y^2&=A(u_1,u_2)\cr
2yy^\prime&=B(u_1,u_2).
}
$$
Moreover the $S$-integer $y$ divides both $F(u_1)$ and $G(u_2)$
in the ring $\O_{S\cup T}$.
}

\medskip

\noindent{\it Proof}.
Of course, equation (3.3) is exactly $y^2=A(u_1,u_2)$ and
$B(X,Y)=\uno X{\partial \over\partial X}A(X,Y)+
\due Y{\partial \over\partial Y}B(X,Y)$, so by the above lemma
we have $2yy^\prime=B(u_1,u_2)$ as wanted.

To prove that $y$ divides
$F(u_1)$ (resp. $G(u_2)$) we show that $F(X)$ (resp. $G(Y)$)
are linear combinations of $A(X,Y),B(X,Y)$ with coefficients in
$\O_{S\cup T}[Y]$ (resp. $\O_{S\cup T}[X])$.
This is a general fact in the theory of resultants; in this case
the linear combination $\due A(X,Y)-B(X,Y)$ equals $F(X)$,
proving the first claim. The polynomial $G(Y)$ is obtained
as the linear combination
$$
\eqalign{
&2\uno\left( 2u_2\left(2\uno-\due\right)+(4-\lambda^2)\uno-\lambda u_1\right)
\cdot A(X,Y)\cr
&+\left(2\lambda^2\uno- 2u_2\left(2\uno-\due\right)
+(4-\lambda^2)\uno-\lambda u_1\right)\cdot B(X,Y),
}
$$
concluding the proof.\cvd

Our next goal is to factor the quadratic polynomials $F(X),G(X)$ in a suitable
finite extension of the function field $\kappa(\cc)$; this  function field
extension will be of the form $\kappa(\dd)$ for a cover $\dd\rightarrow\cc$.
The next Lemma will estimate the Euler characteristic of $\dd$.
\medskip

\noindent{\bf Lemma 3.9}. {\it Suppose the polynomials $F(X),G(X)$
defined above are both non-constant.
There exists a cover $\dd\rightarrow\cc$ of degree $\leq 4$ with the following property:
let $U\subset \dd$ be the set formed by the
pre-images of the zeros of the leading and constant coefficients of $F$ and $G$ and
the pre-images of $S$ and $T$. Then the Euler characteristic of $\dd\setminus U$ verifies
$$
\chi_U(\dd)\leq 30\cdot \chi_S(\cc)+5\cdot\max\{0,2g(\cc)-2\}.
$$
Moreover the polynomials $F(X),G(X)$ split in linear factors in the ring $\kappa(\dd)[X]$. }
\medskip

\noindent{\it Proof.}  Let us define the cover $p:\dd\rightarrow\cc$ (where the curve $\dd$
is complete and non-singular) by the property that the field extension
$\kappa(\dd)/p^*(\kappa(\cc))$ be the splitting field of
$F(X)\cdot G(X)$ over $\kappa(\cc)$. Clearly the field $\kappa(\dd)$ is generated over
$p^*(\kappa(\cc))$ by the square roots of the discriminants of $F(X)$ and $G(X)$,
hence has degree $\leq 4$, proving our first contention. The cover $p:\dd\rightarrow\cc$
can be ramified only over the zeros and poles of such discriminants, and the ramification
index is in any point at most two. The poles are contained in $S\cup T$ and the number
of zeros of the discriminants is bounded by their heights.
The  discriminant of  $F(X)$ is
$$
(\lambda^2-4)\left(\due\right)^2-2\lambda^2\due+4\uno\due;
$$
we bound its height by estimating its possible poles, which are either
poles of $u_2^\prime/u_2$ or poles of $u_1^\prime/u_1$. The sum (with multiplicity)
is then bounded by $2H_{\cc}(u_2^\prime/u_2)+H_{\cc}(u_1^\prime/u_1)$.
By Lemma 3.5 this quantity is bounded by $3\chi_S(\cc)$.
The same arguments show that the height of the second discriminant is at most
$4H_{\cc}(u_1^\prime/u_1)+2H_{\cc}(u_2^\prime/u_2)$,  which by Lemma 3.5 is bounded by
$6\chi_S(\cc)$.
Hence the number of ramification points is at most $\sharp(T\cup S)+9\chi_S(\cc).$
The Riemann-Hurwitz genus formula gives
$$
2g(\dd)-2=(\deg p)(2g(\cc)-2)+ \sum_{P\in \dd}(e_P-1)\eqno(3.8)
$$
where $e_P\in\{1,2\}$ is the ramification index of the morphism $p$ at the point $P$.
By the above estimate we have
$$
\sum_{P\in \dd}(e_P-1)\leq \sharp(T\cup S)+9\chi_S(\cc).\eqno(3.9)
$$
Let $U\subset\dd$ be the finite set appearing in the statement and let
$p(U)\subset \cc$ be its image.
Then $\sharp U\leq [\kappa(\dd):p^*(\kappa(\cc))]\cdot \sharp p(U)$;
this inequality, combined with the inequality (3.9) and the equality (3.8) gives
$$
\eqalign{
2g(\dd)-2+\sharp(U)&\leq (\deg p)(2g(\cc)-2+\sharp p(U))+
\sharp(T\cup S)+9\chi_S(\cc)\cr
&=(\deg p)(2g(\cc)-2+\sharp (S\cup T))
+(\deg p) \sharp(p(U)\setminus(S\cup T))\cr
&+\sharp(T\cup S)+9\chi_S(\cc)\cr
&\leq 4 (2g(\cc)-2+\sharp (S\cup T))+4\sharp(p(U)\setminus (S\cup T))
+\sharp(T\cup S)+9\chi_S(\cc),
}
$$
where we have used the non-negativity of the Euler characteristic $2g(\cc)-2+\sharp(S\cup T)
=\chi_{S\cup T}(\cc)$. Observe now that the points in $p(U)$ that are not in $S\cup T$ are just the zeros
of the leading and constant terms in $F(X)$ and $G(X)$, so they are zeros of $u_1^\prime/u_1$,
$u_2^\prime/u_2$, $u_1^\prime/u_1-u_2^\prime/u_2$ or of $2 u_1^\prime/u_1 -u_2^\prime/u_2$;
since all of these
rational functions are of the form $u^\prime/u$ for suitable $S$-units $u$, their heights are
bounded by $\chi_S(\cc)$ (Lemma 3.5); therefore there are at most
$4\chi_S(\cc)$ such points.  From
the above displayed inequality we then obtain
$$
\eqalign{
\chi_U(\dd)&\leq 4\chi_{S\cup T}(\cc)+16\chi_S(\cc)+\sharp(T)+\sharp(S)+9\chi_S(\cc)
\cr
&\leq 29\chi_S(\cc)+5\sharp(T)+\sharp(S)\cr
&\leq 30\cdot \chi_S(\c)+5\max\{0,2g(\cc)-2\},
}
$$
finishing the proof.\cvd

\medskip

We are now able to prove the following

\noindent{\bf Lemma 3.10}. {\it Let $(u_1,u_2,y)$ be a solution of (3.3)
such that $\uno\neq\due,\quad 2\uno\neq\due$ and $u_1,u_2$ are both non
constant.
Let $\dd,U$ be as in Lemma 3.9.
There exist $U$-units $a,b\in\kappa(\dd)$ such that
$$
\left|\max\{H_{\dd}(a), H_{\dd}(b)\}-\max\{H_{\dd}(u_1),H_{\dd}(u_2)\}\right|
\leq 16\cdot\chi_S(\cc)\eqno(3.12)
$$
and
$$
\sum_{v\in \dd\setminus U}\min\{v(a-1),v(b-1)\}\geq {1\over 4}
\cdot \sum_{v\in\dd\setminus  U} v(y).\eqno(3.13)
$$
Moreover, $a=u_1\alpha^{-1}, b=u_2\beta^{-1}$ for suitable
$\alpha,\beta$ with $F(\alpha)=G(\beta)=0$.
}

\smallskip

In the  right-term in (3.13), the rational function $y\in\kappa(\cc)\subset\kappa(\dd)$ is viewed
as a  rational function on $\dd$.
\medskip

\noindent{\it Proof}. Write the polynomials $F(X), G(X)$ in (3.7) as
$$
\eqalign{
F(X)&=\left(\due-\uno\right)(X-\alpha)\cdot (X-\bar{\alpha}),\cr
G(X)&=2\left(2\uno-\due\right)(X-\beta)\cdot(X-\bar{\beta})
}
$$
where $\alpha,\bar{\alpha}$ (resp. $\beta,\bar{\beta}$) are the roots of $F(X)$ (resp.
$G(X)$). Recall now that by Lemma 3.8, $F(u_1)$ and $G(u_2)$ are both divisible by $y$ in the ring
of  $U$-integers and that the leading and constant coefficients of $F(X)$ and $G(X)$ are
$U$-units, so $\alpha,\bar{\alpha},\beta,\bar{\beta}$ too are $U$-units.
Dividing $F(u_1)$ by the $U$-unit $\alpha\bar{\alpha}(\due-\uno)$
and $G(u_2)$ by the $U$-unit $\beta\bar{\beta} 2(2\uno-\due)$ we obtain that
$y$ divides the resulting quotients in the ring of $U$-integers, so
$$
\sum_{v\in\dd\setminus U}\min\{v(u_1\alpha^{-1}-1)+v(u_1{\bar{\alpha}}^{-1}-1),
v(u_2\beta^{-1}-1)+v(u_2{\bar{\beta}}^{-1}-1)\}\geq \sum_{v\in\dd\setminus U}v(y).
$$
Since
$$
\eqalign{
&\sum_{v\in\dd\setminus U}\min\{v(u_1\alpha^{-1}-1)+v(u_1{\bar{\alpha}}^{-1}-1),
v(u_2\beta^{-1}-1)+v(u_2{\bar{\beta}}^{-1}-1)\}\cr
&\leq
\sum_{v\in\dd\setminus U}
\min\{v(u_1\alpha^{-1}-1),v(u_2\beta^{-1}-1)+v(u_2{\bar{\beta}}^{-1}-1)\}\cr
&+\sum_{v\in\dd\setminus U}
\min\{v(u_2{\bar{\alpha}}^{-1}-1),v(u_2\beta^{-1}-1)+v(u_2\bar{\beta}^{-1})\}
\cr
&\leq
\sum_{v\in\dd\setminus U}\min\{v(u_1\alpha^{-1}-1),v(u_2\beta^{-1}-1)\}
+\sum_{v\in\dd\setminus U}\min\{v(u_1\alpha^{-1}-1),v(u_2\bar{\beta}^{-1}-1)\} \cr
&+\sum_{v\in\dd\setminus U}\min\{v(u_1\bar{\alpha}^{-1}-1),v(u_2\beta^{-1}-1)\}
+\sum_{v\in\dd\setminus U}\min\{v(u_1\bar{\alpha}^{-1}-1),v(u_2\bar{\beta}^{-1}-1)\}
}
$$
we obtain that there exist  $U$-units $a\in\{u_1\alpha^{-1},u_2\bar{\alpha}^{-1}\}$
and $b\in\{u_2\beta^{-1},u_2\bar{\beta}^{-1}\}$ such that (3.13) holds.
We now estimate the difference in (3.12). Clearly it is bounded by $
\max\{H_\dd(\alpha),H_\dd(\bar{\alpha}),H_\dd(\beta),H_\dd(\bar{\beta})\}$.
Since $\alpha, \bar{\alpha}$  are roots of $F(X)$, which is defined explicitely in (3.7),
their poles must be either zeros of $\due-\uno$ or be poles of $\due$.  Since the $\cc$-heights
of $\due-\uno$ and of $\due$ is bounded by $\chi_S(\dd)$, their $\dd$-height is $\leq
4\chi_S(\cc)$.  Hence
$$
\max\{H_\dd(\alpha),H_\dd(\bar{\alpha})\}\leq 8 \chi_S(\cc).
$$
For the same reason
$$
\max\{H_\dd(\beta),H_\dd(\bar{\beta})\}\leq 16 \chi_S(\cc).\eqno(3.14)
$$
\cvd
\medskip

We now bound from below the height of $y$ in terms of the heights
of $u_1,u_2$ (so also in terms of $H(a),H(b)$).
To this end we shall use the following lemma, easily deduced from Theorem 1
 in [Z1]:
\medskip

\noindent{\bf Lemma 3.11}. {\it Let $\dd,U$ be as before, $m\geq 2$ be an integer,
$\theta_1,\ldots,\theta_m$ be $U$-units such that no sum $\sum_{i\in I}\theta_i$
vanishes for any nonempty subset $I\subset\{1,\ldots,m\}$.
Then the $U$-integer $\theta_1+\ldots+\theta_m$ satisfies
$$
\sum_{v\in\dd\setminus U} v(\theta_1+\ldots+\theta_m)\geq H_\dd(\theta_1:\ldots:\theta_m)
-{m\choose 2}\chi_U(\dd).
$$
In particular the height $H_\dd(\theta_1+\ldots+\theta_m)$ is bounded from below by
the right side term above.}
\cvd
\medskip

We shall apply the above Lemma with $m=4$ to $y^2$ which is expressed in (3.3)
as a sum of four $U$-units. We obtain, in the case where no subsum in
the expression $u_1^2+\lambda u_1+u_2+1$ vanishes, the following:
\medskip

\noindent{\bf Lemma 3.12}. {\it For every solution $(y,u_1,u_2)$ in $\O_U\times (\O_U^*)^2$
to equation (3.3), with no vanishing subsum on the right-side term,  we have the lower bound
$$
H_\dd(y)\geq \sum_{v\in\dd\setminus U} v(y)\geq \max\{H_\dd(u_1), H_\dd(u_2)\}-6\chi_U(\dd).
$$
}
\smallskip

\noindent{\it Proof}. It suffices to observe that the projective height
$H_\dd(u_1^2:\lambda u_1:u_2:1)$ is
$\geq\max\{2H_{\dd}(u_1), H_{\dd}(u_2)\}$.
\cvd
\medskip

From inequalities (3.13), (3.12) of Lemma 3.10 and from the above Lemma 3.12 we obtain
that for each solution $(y,u_1,u_2)$ of (3.3) there exists a pair of $U$-units $(a,b)$
such that
$$
\sum_{v\in\dd\setminus U}\min\{v(a-1),v(b-1)\}\geq {1\over 4}
\max\{H_\dd(a),H_\dd(b)\}-6 \chi_U(\dd) - 4 \chi_S(\cc).
$$
Since $\chi_S(\cc)\leq \chi_U(\dd)$ we can deduce from the last inequality
$$
\sum_{v\in\dd\setminus U}\min\{v(a-1),v(b-1)\}\geq {1\over 4}
\max\{H_{\dd}(a),H_{\dd}(b)\}-10\chi_U(\dd).\eqno(3.15)
$$

We are now able to apply the Corollary  2.3 from \S 2,
to deduce the following
\medskip

\noindent{\bf Proposition 3.13}. {\it Let  $(y,u_1,u_2)\in\O_S\times(\O_S^\times)^2$
be a solution of equation (3.3) with no vanishing subsum in the right-side term.
Suppose also that $u_1,u_2$ are both non constant and  $\uno\neq \due$, $2\uno\neq\due$.
Let $U,\dd$ be defined in Lemma 3.9 and  $a,b$ be the $U$-units in $\kappa(\dd)$ defined
 in Lemma 3.10. Then either
$$
\max\{H_\cc(u_1),H_\cc(u_2)\}\leq 2^{14}\cdot 35\cdot \chi_S(\cc) \eqno(3.16)
$$
or $a,b$ verify a multiplicative dependence relation of the form
$$
a^r\cdot b^s=1
$$
for integers $(r,s)\in\Z^2 \setminus\{0\}$ with
$$
\max\{|r|,|s|\}\leq 5. \eqno(3.17)
$$
}

\medskip

\noindent{\it Proof}. We suppose (3.16) does not hold, and want to prove that
$a,b$ verify such a multiplicative dependence relation;
since $H_\cc\leq H_\dd$, we have a fortiori the lower bound
$$
\max\{H_\dd(u_1),H_\dd(u_2)\}\geq 2^{14}\cdot 35\cdot \chi_S(\cc)\geq 2^{14} \chi_U(\dd)
$$
the last inequality following easily from Lemma 3.9.
By (3.12) we have
$$
\max\{H_\dd(a),H_\dd(b)\}\geq \max\{H_\dd(u_1),H_\dd(u_2)\}-16\chi_S(\cc)\geq
\max\{H_\dd(u_1),H_\dd(u_2)\}-16\chi_U(\dd)
$$
so we get from the last two displayed inequalities that
$$
\max\{H_\dd(a),H_\dd(b)\}\geq (2^{14}-16)\chi_U(\dd). \eqno(3.18)
$$
Put $H:=\max\{H_\dd(a),H_\dd(b)\}$ and $\chi:=\chi_U(\dd)$; we check that from (3.15) and (3.18)
it follows that
$$
\sum_{v\in\dd\setminus  U}\min\{v(1-a),v(1-b)\} >3\cdot 2^{1/3} H^{2/3}\chi^{1/3}.\eqno(3.19)
$$
In fact, after (3.15) it suffices to show that ${1\over 4} H-10\chi>3\cdot 2^{1/3}
H^{2/3}\chi^{1/3}$ i.e.
$$
H^{2/3}\left({1\over 4}H^{1/3}-3\cdot 2^{1/3}\chi^{1/3}\right)>10\chi.
$$
From (3.18) it follows that ${1\over 4}H^{1/3}-3\cdot 2^{1/3}\chi^{1/3}>3\cdot 2^{1/3}\chi^{1/3}$;
also it  follows that $H^{2/3}>3^{-1}\cdot 2^{-1/3}\cdot 10\chi^{2/3}$, concluding the verification
of (3.19). Then Corollary 2.3 (i) implies that $a,b$ are multiplicative dependent, i.e.
a relation of the form $a^rb^s=1$ holds, for suitable integers $r,s$
not both zero.
Now the part (ii) of the same corollary gives the bound
$$
\sum_{v\in \dd\setminus U}\min\{v(a-1),v(b-1)\}\leq {\max\{H_\dd(a),H_\dd(b)\}
\over \max\{|r|,|s|\}}.
$$
The above bound and the inequality (3.15) and (3.18) give, with the previous notation for $H,\chi$,
$$
{H\over \max\{|r|,|s|\}}>{1\over 4}H-10\chi > {1\over 5} H,
$$
so $\max\{|r|,|s|\}\leq 5$.
\cvd

\medskip

Proposition 3.13 does not give the full conclusion  of Theorem 3.4, since it provides a
multiplicative dependence relation between $a$ and $b$, not $u_1,u_2$.  To end the proof of
Theorem 3.4, we need a last lemma: \medskip

\noindent{\bf Lemma 3.14}. {\it Let $\dd$ be  as before;
let $u_1,u_2\in\kappa(\dd)$ be non zero rational functions.
Let $A(X,Y)\in \kappa[X,Y]$ be a polynomial without multiple factors and let
$B(X,Y)\in\kappa(\dd)[X,Y]$ be defined, in terms of $A(X,Y)$,
$u_1,u_2$, as
$$
B(X,Y)=\uno X{\partial\over\partial X}A(X,Y)
+\due Y{\partial\over\partial Y}A(X,Y).
$$
Let $\alpha,\beta\in\kappa(\dd)^*$ such that $A(\alpha,\beta)=B(\alpha,\beta)=0$. If $u_1\alpha^{-1}$
and $u_2\beta^{-1}$ satisfy a multiplicative dependence relation of the form
$$
\left({u_1\over \alpha}\right)^r\cdot \left({u_2\over\beta}\right)^s=\mu,\quad r,s\in\Z,\quad \mu\in\kappa^*,
\eqno(3.20)
$$
 then either one between $u_1/\alpha$ and $u_2/\beta$ is
constant or $u_1,u_2$ satisfy a multiplicative dependence relation  of the same type, i.e.
$u_1^r\cdot u_2^s\in\kappa$. }
\medskip

\noindent{\it Proof}. If $u_1$ or $u_2$ is constant, the conclusion
follows easily (observing that $\alpha,\beta$ must be constant, since $A(X,Y)$ has no multiple factors); so we suppose them to be non-constant. We also suppose, as we may, that $r,s$ are not both zero.

We first observe that multiplicative dependence relations of the above type correspond to linear relations over the rationals of the form
$$
r\left(\uno-{\alpha^\prime\over\alpha}\right)+s\left(
\due-{\beta^\prime\over\beta}\right)=0.\eqno(3.21)
$$
Starting from the
relation $A(\alpha,\beta)=0$ we obtain, by taking differentials,
$$
\alpha^\prime{\partial\over\partial X}A(\alpha,\beta)+
\beta^\prime{\partial\over\partial Y}A(\alpha,\beta)=0.
$$
Since $B(X,Y)$ too vanishes at the point $(\alpha,\beta)$ we get,
applying the definition of $B(X,Y)$,
$$
\uno\alpha  {\partial\over\partial X}A(\alpha,\beta)+
\due\beta {\partial\over\partial Y}A(\alpha,\beta)=0.
$$

Suppose first that $(\alpha,\beta)$ is singular on the curve $A(X,Y)=0$. Then $\alpha,\beta$ are both  in $\kappa$, because $A(X,Y)$ has coefficients in $\kappa$ and  no multiple factors;  in this case the conclusion of the lemma is automatic.
Otherwise,
 the above two relations 
 imply that the matrix $\left(\matrix{\uno&\due\cr {\alpha^\prime\over
\alpha}&{\beta^\prime\over\beta}}\right)$
has rank one, so the same holds for the matrix
$$
 \left(\matrix{
\uno-{\alpha^\prime\over \alpha}&\due-{\beta^\prime\over\beta}\cr
{\uno}&{\due}
} \right).
$$
If the first row in this matrix vanishes we are done, since then $u_1/\alpha$ (and $u_2/\beta$) would be constant, falling into one alternative of the conclusion of the lemma. Then we shall suppose that the first row does not vanishes. Since it is
orthogonal to the vector $(r,s)$ by (3.21), and since the matrix has rank one, the same vector $(r,s)$ must be orthogonal also to the second row,  obtaining the sought conclusion that $r\uno+s\due=0$.\cvd

\medskip

\noindent{\bf Proof of Theorem 3.4}.  We shall suppose that conditions
(i), (ii) and (iii) of Theorem 3.4 are not satisfied and deduce
a contradiction. Since (i) is not  satisfied no
subsum in the right term of (3.3) can vanish.
Also, the polynomial $F(X)$ defined by (3.7) has
degree exactly two, in particular is not constant, because the vanishing of its leading
coefficient (which is $\uno-\due$) would imply a multiplicative relation of the form
$u_1=\mu\cdot u_2$, for a constant $\mu\in\kappa^*$, which has been excluded by (ii).
The same remarks apply to the leading coefficient of $G(X)$, which vanishes
only if $u_1\cdot u_2^{-2}$ is a constant. Also $\uno\neq 0$ and $\due\neq 0$, otherwise
$u_1$ or $u_2$ would be constant, hence the pair $(u_1,u_2)$ would satisfy a multiplicative
dependence relation with exponents $(1,0)$ or $(0,1)$.
Then neither $F(X)$ nor $G(X)$ are constant polynomials and their constant
coefficient are nonzero.

Then Proposition 3.13 and 
 the above lemmas apply and we obtain that for a suitable root
$\alpha$, (resp. $\beta$) of the polynomial $F(X)$ (resp. $G(X)$)
defined in (3.7), the $U$-units $a=u_1\alpha^{-1}$ and $b=u_2\beta^{-1}$
satisfy a multiplicative dependence relation as in Proposition 3.13.
\medskip

Now our formulas (3.7) yield $F(X)=c_1A(X,Y)-B(X,Y), G(Y)=c_2A(X,Y)+c_3B(X,Y)$ for certain $c_1,c_2,c_3\in\kappa(\cc)$. Now, we have $c_1c_3+c_2\neq 0$, for otherwise e$F(X), G(Y)$ would be linearly dependent over $\kappa(\cc)$, and thus neither depend on $X$ nor $Y$, which is impossible (looking at their coefficients, we would find that $u_1,u_2$ are both constant). Then the vanishing of $F(\alpha), G(\beta)$ implies $A(\alpha,\beta)=B(\alpha,\beta)=0$.

But then, by the above Lemma, either $u_1/\alpha$ is constant, or $u_2/\beta$ is constant
or $u_1,u_2$ too satisfy a multiplicative  dependence relation (with the same exponents as
for $a$ and $b$). The first two possibilities are ruled out by height consideration:
in fact if $u_1/\alpha$ were constant, the height of $u_1$ would be equal to the height of
$\alpha$ which is bounded by the inequality (3.14);
 since we assume that (iii) of Theorem 3.4 does not hold, this is excluded. For
the same reason $u_2/\beta$ is not constant. The last possibility coincides with (ii),
which we have excluded, concluding the proof.  \cvd

\medskip

We now prove Theorem 3.1. which, as we already remarked, is equivalent to Theorem 1.1.
\medskip

\noindent{\bf Proof of Theorem 3.1}. By Lemma 3.3 we can
find a scalar $\lambda\in\kappa$ and identify $X$ with the open subset of the plane $\A^2$
defined by
$$
x\cdot (y^2-x^2-\lambda x-1)\neq 0.
$$
Then, as in Lemma 3.3, the morphism $f:\cc\setminus S\rightarrow X$
 will be expressed in the form
$$
f:p\mapsto (u_1(p),y(p))
$$
where $(y,u_1)\in\O_S\times(\O_S)^*$ satisfies (3.3) for a suitable
$S$-unit $u_2\in\O_S^*$.
Then one of the conclusions (i),...,(iii) of Theorem 3.4 holds.
We shall prove the estimate of Theorem 3.1 in each case:
\smallskip

\item{(i)} suppose that some subsum on the right term of (3.3) vanishes;
 then the curve $f(\cc)$ is either a line or a conic, hence its degree is $\leq 2$.
\smallskip

\item{(ii)}if $u_2^s=u_1^r\cdot \mu$, where $\mu\in\kappa$ and
the pair of integers $(r,s)\in\Z^2\setminus\{0\}$ satisfy
$\max\{|r|,|s|\}\leq 5$, then the degree of $f(\cc)$ is $\leq 20$;
\smallskip

\item{(iii)} finally, if $\max\{H_{\cc}(u_1),H_{\cc}(u_2)\}\leq 2^{14}\cdot 35\cdot\chi(\c)$,
then by Lemma 3.2 the estimate of Theorem 3.1 holds.
\cvd

\bigskip

\noindent{\bf Proof of Theorem 1.2}.
We now prove the full Theorem 1.2, of which Theorem 3.4 (hence
also Theorems 1.1 and 3.1) is just a particular case. We avoid however the rather lengthy explicit
calculations of the constants involved in the estimates.
\smallskip

We begin by factoring the polynomial $A(X,Y)\in\kappa[X,Y]$ as 
$$
A(X,Y)=A_1(X,Y)\cdots A_l(X,Y),
$$
into irreducible factors $A_1(X,Y),\ldots, A_l(X,Y)\in\kappa[X,Y]$. Now, given a pair of $S$-units
$u_1,u_2$, every multiple zero (outside $S$) of the rational function $A(u_1,u_2)\in\kappa(\cc)$  
is either a common zero to two factors $A_i(u_1,u_2),A_j(u_1,u_2)$, for a pair $1\leq i<j\leq l$,
or is a multiple zero of one of the factors. 
Hence Theorem 1.2 is a consequence of the following two assertions,  of which the second one is
just  Theorem 1.2 for irreducible polynomials:

\medskip

\noindent {\bf Proposition 3.15}. {\it Let $A_1(X,Y),A_2(X,Y)\in\kappa[X,Y]$ be coprime polynomials,
$\epsilon>0$ a positive real number. There exists an integer $C_2=C_2(A_1,A_2,\epsilon)$ such
that the solutions  $(u_1,u_2)\in\O_S^*$ to the inequality      
$$
\sum_{v \in\cc\setminus S}\min\{v(A_1(u_1,u_2)),v(A_2(u_1,u_2))\} > \epsilon\max\{H(u_1),H(u_2)\}
\eqno(3.22)  
$$
either verify $\max\{H(u_1),H(u_2)\}\leq C_2\cdot\max\{1,\chi(\cc\setminus S)\}$ or verify a
multiplicative dependence relation of the form $u_1^r u_2^s\in\kappa$, for some nonzero vector
$(r,s)\in\Z^2$.
} \medskip

Note that the particular case when $A_1(X,Y)=X-1$, $A_2(X,Y)=Y-1$
consists exactly of Corollary 2.3. \medskip

\noindent {\bf Proposition 3.16}. {\it Let $A(X,Y)\in\kappa[X,Y]$ be an irreducible
polynomial.   For every positive $\epsilon$ there exists an integer
$C_3=C_3(A,\epsilon)$ such that 
for all pairs $(u_1,u_2)\in(\O_S^*)^2$ with $\max\{H(u_1),H(u_2)\}\geq
C_3\cdot\max\{1,\chi(\cc\setminus S)\}$, either the rational function $A(u_1,u_2)\in\kappa(\cc)$
satisfies 
$$
\sum_{v\in\cc\setminus S} \max\{0, v(A(u_1,u_2))-1\}<\epsilon \max\{H(u_1),H(u_2)\},\eqno(3.23)
$$
or $u_1,u_2$ verify a multiplicative dependence relation of the form
$$
u_1^r\cdot  u_2^s\in\kappa 
$$
for a suitable pair of integers $(r,s)\in\Z^2\setminus\{0\}$.
}
\medskip

\noindent{\it Proof of Proposition 3.15}.  Let $F(X)$ (resp.
$G(Y)$) be the resultant of the polynomials $A_1(X,Y),A_2(X,Y)$ with respect to $Y$ (resp. $X$);
they are both non-zero since by assumption $A_1(X,Y),A_2(X,Y)$ are coprime. Since $F(X),G(Y)$ are
linear combinations of $A_1(X,Y),A_2(X,Y)$ with coefficients in $\O_S[X,Y]$, every common zero of
$A(u_1,u_2),A_2(u_1,u_2)$ outside $S$ is also a common zero of $F(u_1),G(u_2)$, with at least the
same multiplicity.  This means that the left-hand side of (3.22) verifies   
$$
\sum_{v\in \cc\setminus S}\min\{v(A_1(u_1,u_2)),v(A_2(u_1,u_2))\}\leq\sum_{v\in\cc\setminus S}
\min\{v(F(u_1)),v(G(u_2))\}.
$$
As in Lemma 3.9, we can factor the polynomials $F(X)$, $G(Y)$, but in this case both $F(X),
G(Y)$ have constant coefficients, so they factor already in the rings $\kappa[X],\kappa[Y]$
respectively; we then write 
$$
F(X)=\alpha\cdot X^e\cdot (X-\alpha_1)\cdots(X-\alpha_m),\qquad G(Y)=\beta\cdot Y^f\cdot
(Y-\beta_1)\cdots(Y-\beta_n), $$
for suitable $a,b\in\kappa^*,\alpha_1,\ldots,\alpha_m,\beta_1,\ldots,\beta_n\in \kappa^*$ and
non-negative integers $e,f$.  As already remarked, 
$$
\sum_{v\in\cc\setminus S}
\min\{v(F(u_1)),v(G(u_2))\}\leq  m  n \cdot \max_{(i,j)}\sum_{v\in \cc\setminus
S}\min\{v(u_1-\alpha_i),v(u_2-\beta_j)\},
$$
where the maximum is taken over the pairs $(i,j)\in\{1,\ldots, m\}\times\{1,\ldots, n\}$.
We then fix a pair $(i,j)$ and estimate the corresponding quantity $\sum_{v\in \cc\setminus
S}\min\{v(u_1-\alpha_i),v(u_2-\beta_j)\}$. Corollary 2.3, applied with $a=u_1/\alpha_i$,
$b=u_2/\beta_j$ immediately implies that either $a,b$ are multiplicative dependent modulo
constants (so are $u_1,u_2$), or the bound
$$
\sum_{v\in \cc\setminus
S}\min\{v(u_1-\alpha_i),v(u_2-\beta_j)\}\leq
3\root{3}\of{2}\max\{H(u_1),H(u_2)\}^{2/3}\chi(\cc\setminus S)^{1/3}
 $$
holds for all $(i,j)$. In the first case we have our content; suppose now the above inequality
holds: then the left-hand side of (3.22) will be bounded as
 $$
\sum_{v \in\cc\setminus S}\min\{v(A_1(u_1,u_2)),v(A_2(u_1,u_2))\} \leq mn
\cdot 3\root{3}\of{2}\max\{H(u_1),H(u_2)\}^{2/3}\chi^{1/3}.
 $$
The right-hand side term will be $<\epsilon \max\{H(u_1),H(u_2)\}$ provided  
$ \max\{H(u_1),H(u_2)\}>mn\cdot 3^3\cdot 2\cdot\epsilon^{-3}\chi$, so we have proved the Proposition with $C_2=
 mn\cdot 3^3\cdot 2\cdot\epsilon^{-3}$. Note that it is indeed a function of
$A_1(X,Y),A_2(X,Y)$ and $\epsilon$ only.
\cvd 
\medskip

\noindent{\it Proof of Proposition 3.16}.  
Suppose $A(X,Y)$ is irreducible, and that $A(u_1,u_2)$ has ``many" multiple zeros,  with
$S$-units $u_1,u_2\in\O_S^*$. Taking the differential of $A(u_1,u_2)$, and using again the above
formalism for the operator ${}^\prime$, we obtain a second function $B(u_1,u_2)$ with many zeros in
common with $A(u_1,u_2)$, where the polynomial $B(X,Y)\in\kappa(\cc)[X,Y]$ is defined in
the statement of Lemma 3.7, namely 
$$
B(X,Y)=\uno X{\partial A\over\partial X}+\due Y{\partial A\over\partial Y}.
$$
We enlarge the set $S$ to a set $U_1$ including the set $T$ of Lemma 3.5, so that the coefficients of $B$ have no poles outside $U$. Note that the enlargement is harmless for our purposes, since the set $T$ contains at most $\max\{1,\chi(\cc\setminus S)\}$ points. So $\chi(\cc\setminus U_1)\leq 2\chi(\cc\setminus S)$.

We have, since $B(u_1,u_2)=A(u_1,u_2)^\prime$,
$$
\sum_{v\in\cc\setminus U_1}\max\{0,v(A(u_1,u_2))-1\}\leq \sum_{v\in\cc\setminus U_1}
\min\{v(A(u_1,u_2)),v(B(u_1,u_2))\}, \eqno(3.24) 
$$
so it suffices to bound the right-hand side of (3.24), in order to prove Proposition 3.16.
To do this we subdivide the proof into six steps.
 \medskip

{\bf Step 1}. {\it If the $S$-units $u_1,u_2$ are multiplicatively independent modulo $\kappa^*$, then the polynomials $A(X,Y)\in\kappa[X,Y]\subset\kappa(\cc)[X,Y]$  and 
$B(X,Y)\in\kappa(\cc)[X,Y]$ are coprime.} \medskip

In fact, since $A(X,Y)$ is irreducible and $\deg
B(X,Y)\leq\deg A(X,Y)$, if they were not coprime there would exist an element
$\lambda\in\kappa(\cc)$ such that $B(X,Y)=\lambda A(X,Y)$. Now observe that the operator $\uno
X{\partial\over\partial X}+\due Y{\partial\over\partial Y}$ sends the monomial $X^iY^j$ to the
monomial $(i\uno +j\due)X^iY^j$. Hence from $B(X,Y)=\lambda A(X,Y)$ it would follow that for every
$(i,j)$ such that the monomial $X^iY^j$ appears in $A(X,Y)$, $i\uno+j\due=\lambda$; considering two
different monomials\note{if $A(X,Y)$ were a monomial, the function $A(u_1,u_2)$ would have no
zero outside $S$, and (3.23)  would be trivial} we obtain
$i\uno+j\due=h\uno+k\due$, for a choice of pairs $(i,j)\neq(h,k)$, hence $(i-h)\uno=(k-j)\due$,
which implies the multiplicative dependence of $u_1,u_2$.

Now, taking again the resultants $F(X)={\rm Res}_Y(A(X,Y),B(X,Y))$, $G(Y)={\rm
Res}_X(A(X,Y),B(X,Y))$, which do not vanish identically in view of the coprimality of $A(X,Y), B(X,Y)$,  we obtain an analogue of the last conclusion of Lemma 3.8, so a ``large" $\gcd$ between $F(u_1)$ and $G(u_2)$. 
More precisely, since $F(X), G(Y)$ are linear combinations of $A(X,Y), B(X,Y)$ over $\O_{U_1}[X,Y]$, for every place $v\not\in U_1$ 
  $$
\min\{v(A(u_1,u_2)),v(B(u_1,u_2))\}\leq 
\min\{v(F(u_1)),v(G(u_2))\}.\eqno(3.25)
 $$ 

 So, in view of (3.24), we are reduced to bound the $\gcd$ of $F(u_1),G(u_2)$, as in the Proof of Proposition 3.15. The only difference is that now   the
polynomials $F(X),G(Y)$ have their coefficients in the ring $\O_{U_1}\subset\kappa(\cc)$, and that their coefficients also involve the $S$-units $u_1,u_2$.  For this reason we need steps 2 and 3.
\medskip

{\bf Step 2}. {\it With the above notation, there exist  numbers
$C_4=C_4(A)>0, C_5=C_5(A)\geq 1$ such that the height of $F(X),G(Y)$ is bounded by
$C_4\cdot\max\{1,\chi(\cc\setminus S)\}$ and the degrees of $F(X),G(Y)$ are bounded by $C_5$. }
 \smallskip

\noindent{\it Proof}. Let us consider the resultant $F(X)$, with respect to $Y$; the same
reasoning will apply to $G(Y)$. One can express explicitly  $F(X)$ as the determinant of a
$N\times N$ square matrix, where $N=\deg_Y A(X,Y)+\deg_Y B(X,Y)\leq 2\deg_Y A(X,Y)$ only depends on
$A(X,Y)$; the entries of such a matrix are the monomials appearing in $A(X,Y),B(X,Y)$, so their height is bounded in virtue of Lemma 3.5 by $H(\uno)+H(\due)\leq 2\chi(\cc\setminus S)$. Then the height of $F(X)$ is  bounded by $2N\cdot \max\{1,\chi(\cc\setminus S)\}$, and our first contention is proved by taking  
$$
C_4=2N=4\max\{\deg_XA(X,Y), \deg_Y A(X,Y)\}.
$$  
Now, it is clear that the degree of $F(X)$ is $\leq N\cdot \max\{\deg_X
A(X,Y),\deg_X B(X,Y)\}=N\deg_X A(X,Y)$, so the assertion is proved after taking $C_5=2\deg_X
A(X,Y)\deg_Y A(X,Y)$. \cvd

\medskip

{\bf Step 3}. As in the proof of Proposition 3.15 (and of Theorem 1.1), we have to factor
$F(X),G(Y)$ in a suitable finite extension of $\kappa(\cc)$; we need to estimate the genus of their
splitting field and the height of their roots;  the analogue of Lemma 3.9 is the following fact 
\medskip

{\it There exists a number  $C_6=C_6(A)$ with the following property:
There exists a curve $\dd$, a cover $\dd\rightarrow \cc$, of degree $\leq 2C_5$, a finite set
$U_2\subset \dd$ such that the polynomials $F(X),G(Y)$ split over $\kappa(\dd)$ in linear factors;
also, their roots are $U_2$-units and the Euler characteristic of $\dd\setminus U_2$ verifies
$$
\chi(\dd\setminus U_2)\leq C_6\cdot \max\{1,\chi(\cc\setminus S)\} .\eqno(3.26)
$$
}

\smallskip

 The proof of this statement is very similar to the one of Lemma 3.9. Let us
define the cover $p:\dd\rightarrow\cc$ (where the curve $\dd$ is complete and smooth) by the
property that  $\kappa(\dd)/\kappa(\cc)$ be the splitting field of $F(X)G(X)$ over $\kappa(\cc)$.
As already remarked,  the degree of $p$ is $\leq 2C_5$. The ramification can arise only over
zeros or poles of the discriminants of the irreducible factors of $F(X)\cdot G(X)$.  The poles are
contained in $S\cup T$, for the finite set $T$ defined in Lemma 3.5; to bound the number of zeros,
observe that the discriminants of (the irreducible factors  of) $F(X),G(X)$, viewed as polynomials
in the coefficients of $F(X),G(X)$, have degree $<(\deg(F(X))^2$  (resp. $\leq (\deg G(X))^2$),
hence the degree of the ramification divisor (in $\cc$) is bounded by  $2C_5^2+\sharp(S\cup T)\leq 2
C_5^2+\sharp(S)+\max\{0,2g(\cc)-2\}$. Observing that each ramification index is $\leq C_5$ we
obtain that  the total ramification of $p$ (i.e. the sum of the ramification indices) is bounded by
$2C_5^3+C_5\cdot (\sharp(S)+\max\{0,2g(\cc)-2\})$. By the Riemann-Hurwitz formula, the Euler
characteristic of  $\dd$ verifies $$
\eqalign{
\chi(\dd)&\leq \deg(p)\chi(\cc) + 2C_5^3+C_5\cdot (\sharp(S)+\max\{0,2g(\cc)-2\})
\cr &\leq
C_5(\chi(\cc)+\sharp(S)+\max\{0,2g(\cc)-2\})+2C_5^3\leq 2C_5\max\{1,\chi(\cc\setminus S)\}+2C_5^3.
}
$$
Now define $U_2\subset\dd$ to be the pre-image   $p^{-1}(S\cup T\cup U^\prime)$, where $U^\prime$ is the set of zeros of the constant and leading coefficients of $F(X)$ and of $G(X)$. 
The cardinality of $U^\prime$ is bounded by $2H(F(X))+2H(G(X))\leq 4  C_4$; this bound for $\sharp(U^\prime)$ and the above one for $\chi(\dd)$ easily imply
 $$
\chi(\dd\setminus U_2)\leq 3\cdot C_5^3\max\{1,C_4\}\cdot \max\{1,\chi(\cc\setminus S)\}=
C_6\cdot\max\{1,\chi(\cc\setminus S)\},
 $$
where $C_6=C_6(A)$ only depends on the polynomial $A(X,Y)$, as wanted. \cvd
 \medskip

We now factor $F(X)$ (resp. $G(Y)$) in the ring $\kappa(\dd)[X]$ as
$$
F(X)=\lambda\cdot (X-\alpha_1)\cdots (X-\alpha_m),\qquad G(Y)=\mu\cdot (Y-\beta_1)\cdots
(Y-\beta_n), 
$$
for $U_2$-units $\lambda,\mu,\alpha_1,\ldots,\alpha_m,\beta_1,\ldots,\beta_n$ and the degrees $mn$
are $\leq C_5$. 
\medskip

{\bf Step 4}. In order to bound (from above) the right-hand side of (3.25) for $v\in\cc\setminus U_1$, which is our goal, we begin by estimating  
$ \sum_{v\in\dd\setminus U_2}\min\{F(u_1),G(u_2)\}.$

Working on the curve $\dd$, we have for every place $v\in\dd\setminus U_2$,
$$
\min\{v(F(u_1)),v(G(u_2))\}\leq \sum_{(i,j)} \min\{v(u_1-\alpha_i),v(u_2-\beta_j)\} 
$$
(This follows from the inequality $\min\{\sum_i a_i,\sum_j b_j\}\leq \sum_{i,j}\min\{a_i,b_j\}$, holding for arbitrary real non-negative numbers $a_i,b_j$.) 
From $(3.25)$ we have, for every $v\in\dd\setminus U_2$,
$$
\min\{v(A(u_1,u_2)), v(B(u_1,u_2)\}\leq \sum_{(i,j)} \min\{v(u_1-\alpha_i),v(u_2-\beta_j)\}. \eqno(3.27). 
$$

Let us define the set  ${\cal Z}$  by
$$
{\cal Z}=\{(\alpha,\beta) \, :\, A(\alpha,\beta )=B(\alpha ,\beta)=0\}\qquad
$$
We define the set $U$ by enlarging $U_2$  so that the non-zero ones among the values $A(\alpha_i,\beta_j), B(\alpha_i,\beta_j)$ are $U$-units; namely we add to $U_2$ the zeros in $\dd$ of those functions $A(\alpha_i,\beta_j), B(\alpha_i,\beta_j)$ which do not vanish.
 Note that  by the height estimates  carried out in  Step 2 and Step 3 the heights and degrees of $F(X), G(Y)$ are bounded by a constant time $\max\{1,\chi(\cc\setminus S)\}$; therefore the heights of $\alpha,\beta$, $A(\alpha,\beta), B(\alpha,\beta)$ are likewise bounded. In turn, this implies that the number of zeros of nonvanishing elements $A(\alpha_i,\beta_j), B(\alpha_i,\beta_j)$ is at most $\gamma^\prime \max(1,\chi(\cc\setminus S))$, for a suitable constant $\gamma^\prime$  depending only on $m,n$. \medskip

We now show that $(3.27)$ holds even after restricting the summation on the right to the pairs $(\alpha_i,\beta_j)\in{\cal Z}$. To do this we may suppose that 
$\min\{v(A(u_1,u_2)), v(B(u_1,u_2)\}>0$, so $A(u_1,u_2)\equiv B(u_1,u_2)\equiv 0$ $\pmod v$.  Take a point $(\alpha_i,\beta_j)$ outside ${\cal Z}$. 
By definition of ${\cal Z}$, at least one between $A(\alpha_i,\beta_j), B(\alpha_i,\beta_j)$  is non-zero, say $A(\alpha_i,\beta_j)\neq 0$;  by the present choice of $U$, we deduce that $A(\alpha_i,\beta_j)$  is a $U$-unit, so in particular $v(A(\alpha_i,\beta_j))=0$. Then we must have $\min\{v(u_1-\alpha ),v(u_2-\beta)\}=0$, for otherwise we would have $\alpha_i\equiv u_1, \beta_j\equiv u_2$ $\pmod v$, so $A(\alpha_i,\beta_j)\equiv A(u_1,u_2)\equiv 0$ $\pmod v$, contradicting the fact that $A(\alpha_i,\beta_j)$ is a $U$-unit. In conclusion, if the term on the left in $(3.27)$ is positive, the terms on the right vanish outside ${\cal Z}$, so the inequality holds on the restricted summation. 
If the term on the left is zero, the inequality holds trivially on the restricted summation, since every term is non-negative. We write explicitly this conclusion as
$$
\min\{v(A(u_1,u_2)), v(B(u_1,u_2)\}\leq \sum_{(\alpha,\beta)\in{\cal Z}} \min\{v(u_1-\alpha),v(u_2-\beta)\}\qquad {\rm for\ all}\ v\in\dd\setminus U.\eqno(3.27b)
$$
We may now sum over $U$ and invert the order of summation, obtaining
$$
\eqalign{
\sum_{v\in\cc\setminus U}\min\{v(A(u_1,u_2),v(B(u_1,u_2)\}&\leq
\sum_{(\alpha,\beta)\in{\cal Z}}\sum_{v\in\dd\setminus U} \min\{v(u_1-\alpha),v(u_2-\beta)\}\cr
&\leq mn
\max_{(\alpha,\beta)\in {\cal Z}}\sum_{v\in\dd\setminus U}\min\{v(u_1-\alpha ),v(u_2-\beta)\}.} \eqno(3.27c)
$$
\medskip

Fix the pair $(\alpha,\beta)\in{\cal Z}$ attaining the maximum; 
then Corollary 2.3, applied with $a=u_1/\alpha, b=u_2/\beta$ and $(\dd,U)$ instead of $(\cc,S)$, gives the bound
$$
\sum_{v\in\dd\setminus U} \min\{v(u_1-\alpha),v(u_2-\beta)\}\leq
3\root{3}\of{2}\max\left\{H_{\dd}\left({u_1\over
\alpha}\right),H_{\dd}\left({u_2\over\beta}\right)\right\}^{2/3}\cdot  \chi(\dd\setminus
U)^{1/3},\eqno(3.28)
 $$
unless $u_1/\alpha, u_2/\beta$ are multiplicatively dependent modulo constants. In this last case, by Lemma 3.14, either $u_1,u_2$ are also multiplicatively dependent modulo constant, and we are done, or one between $u_1/\alpha, u_2/\beta$ is constant. If this constant is $\neq 1$, the above upper bound holds trivially (the sum on the left side of $(3.28)$ is zero). Otherwise, say that $u_1=\alpha$. Now the polynomial $A(\alpha,Y)$ has no multiple factors and one is reduced to estimate the number of multiple zeros of $u_1-\rho$, for a root $\rho$ of $A(\alpha,Y)=0$.  This may be done by ususal $abc$ and leads to the conclusion of Proposition (3.13).\smallskip

(We point out to the reader that all the previous efforts to restrict the summation over ${\cal Z}$ were necessary merely to apply Lemma 3.14.)
\smallskip

Therefore we assume from now on that $u_1/\alpha, u_2/\beta$ are multiplicatively independent modulo constants.

It remains to show that $(3.28)$  implies $(3.23)$, provided the maximum height  of $u_1,u_2$ is sufficiently large.  First of all, by $(3.26)$, we can rewrite $(3.28)$ as
$$
 \sum_{v\in\dd\setminus U} \min\{v(u_1-\alpha),v(u_2-\beta)\}\leq
C_7\cdot\max\{H_\dd(u_1/\alpha),H_\dd(u_2/\beta)\}^{2/3}\cdot\max\{1,\chi(\cc\setminus
S)\}^{1/3},\eqno(3.29) 
$$
for a suitable $C_7=C_7(A)$.  Note that this essentially implies, via $(3.24), (3.27c)$, the inequality $(3.23)$ of Proposition 3.16. The only difference is represented by the sum running over the restricted set of places outside $U$, so $(3.23)$ is proved for the moment only for the corresponding sum over $\cc\setminus p(U)\subset\cc\setminus S$. 
In the following two steps we deal with the easy task of estimating the contribution in $p(U)\setminus S$ and to pass from $\dd$ to $\cc$. 
\medskip

{\bf Step 5}. We now estimate the sum in $p(U)\setminus S $ of $\max(0,v(A(u_1,u_2)-1)$.
We write $A(u_1,u_2)$ as a sum $\theta_1+\ldots+\theta_M$ of $S$-units, which are monomials in $u_1,u_2$. We may suppose that no subsum vanishes. Note that $H_{\cc}(\theta_1:\ldots:\theta_M)\geq H_\cc(A(u_1,u_2))=\sum_{v\in\cc}\max\{0,v(A(u_1,u_2))\}$. Then by Lemma 3.11, (with $\cc$ in place of $\dd$, $p(U)$ in place of $U$ and $M$ in place of $m$) we deduce that 
$$
\sum_{v\in p(U)}\max\{0,v(A(u_1,u_2))\}\leq {M\choose 2}\chi(\cc\setminus p(U)).
$$

{\bf Step 6}. Finally by $(3.24), (3.27c), (3.29)$ we have
$$
\sum_{v\in\dd\setminus U}\max(0,v(A(u_1,u_2))-1)\leq C_7\cdot\max\{H_\dd(u_1/\alpha),H_\dd(u_2/\beta)\}^{2/3}\cdot\max\{1,\chi(\cc\setminus
S)\}^{1/3}. 
$$
Viewing the places and heights relative to $\cc$ rather then $\dd$ gives 
$$
\sum_{v\in\cc\setminus p(U)}\max(0,v(A(u_1,u_2))-1)\leq C_8\cdot\max\{H_\cc(u_1/\alpha),H_\cc(u_2/\beta)\}^{2/3}\cdot\max\{1,\chi(\cc\setminus
S)\}^{1/3}. 
$$
By  Step 5 
$$
\sum_{v\in\cc\setminus S}\max(0,v(A(u_1,u_2))-1)\leq C_9\cdot\max\{H_\cc(u_1/\alpha),H_\cc(u_2/\beta)\}^{2/3}\cdot\max\{1,\chi(\cc\setminus
S)\}^{1/3} .
$$
This concludes the proof.

\medskip

\noindent {\bf Proof of Theorem 1.3}.  It can be easily deduced from Theorem 1.2. In fact, the function field of $Y$ is a cyclic extension of the function field of $\Gm^2$, so it is obtained by taking the $d$-th root of a suitable rational function on $\Gm^2$. Since by hypothesis there exists a finite map $\pi:Y\rightarrow\Gm^2$ we can write the $\kappa$-algebra $\kappa[Y]$ as the extension of the $\kappa$-algebra $\kappa[\Gm^2]=\kappa[u_1,u_2]$ obtained by adding a $d$-th root:
$$
\kappa[Y]=\kappa\left[u_1,u_2,{1\over u_1},{1\over u_2}\right][\root{d}\of{F(u_1,u_2)}]
$$
for a suitable polynomial $F(X,Y)\in\kappa[X,Y]$. Our hypothesis on the log-Kodaira dimension of
$Y$ (equivalent to the hypothesis that the morphism $\pi:Y\rightarrow\Gm$ is ramified) implies
that  $F(X,Y)$ is not a perfect $d$-th power. Let now $\cc$ be a curve, $S\subset\cc$ a
finite set as before and $f:\cc\setminus S\rightarrow Y$ a morphism. Composing it with $\pi$ we
obtain a morphism $\pi\circ f:\cc\setminus S\rightarrow\Gm^2$ which is expressed as $(\cc\setminus
S)\ni p\mapsto (u_1(p),u_2(p))\in\Gm^2$ for suitable $S$-units $u_1,u_2$ in $\kappa(\cc)$. The fact
that such a morphism factors through $Y$ just means that the rational function $F(u_1,u_2)$ is a
perfect $d$-th power, so in particular all its zeros have multiplicity divisible by $d$. Let us now
factor  $F(X,Y)$ in the ring $\kappa[X,Y]$ as $$
F(X,Y) = F_1(X,Y)^{e_1}\cdots F_l(X,Y)^{e_l}\cdot G(X,Y)^d
$$
where $F_1(X,Y),\ldots,F_l(X,Y)$ are pairwise coprime, $e_1,\ldots, e_l$ are integers all $<d$.   

Since the polynomial $F(X,Y)$ is not a perfect $d$-th power, the factor $A(X,Y):=F_1(X,Y)\cdots
F_l(X,Y)$ in non-constant, and it is easily seen that for every zero $v\in\cc\setminus S$ of
$F(u_1,u_2)$  one has $v(A(u_1,u_2))>1$. Then Theorem 1.2 applies and we deduce either a bound on
the height of such a morphism $f:\cc\setminus S\rightarrow Y$ or a multiplicative dependence
relation between $u_1$ and $u_2$. In both cases one obtains a bound for the  projective degree of
the image $\pi\circ f(\cc)$ in $\Gm^2\subset\Pr_2$, which in turn gives a bound for the projective
degree of the curve $f(\cc)$ with respect to the given  embedding.  
\cvd

\bigskip

\noindent {\bf \S 4. Maps to the complement of a cubic}.
\smallskip

Since the canonical divisor of $\Pr_2$ has degree $-3$, the complement of a smooth cubic
$D\subset \Pr_2$, or of a singular cubic with normal crossing singularities, has  vanishing log
Kodaira dimension. Hence one expects that the conclusion of Theorem 1.1 does not hold in this case, so the degree of a curve in $\Pr_2\setminus D$ cannot be bounded in term of its Euler characteristic. We prove that this is the case.

\medskip

We distinguish three cases, according to the number of components of $D$. 

First consider the case of an irreducible cubic $D$.\smallskip

{\bf Theorem 4.1}. {\it Let $D\subset\Pr_2$ be an irreducible cubic. For every integer $d\geq 1$ there exists a curve $Y\subset X:=\Pr_2\setminus D$ such that $\deg Y\geq d$ and $\chi(Y)\leq 12$.}
\medskip

We start with the easiest cases of singular cubics.

Consider first the case when $D$ has a cusp. Then in suitable coordinates $D$ is given by the equation $D: zy^2-x^3=0$. Given a polynomial $p(t)\in\kappa[t]$, with $p(0)\neq 0$, and a positive integer $n$, the map $\Gm\to \Pr_2$  sending $\Gm\ni t\mapsto (t^{2n}p(t):t^{3n}:p(t)^3+1)$ avoids the divisor $D$ and its image is a curve of degree $\geq 2n+\deg(t)$.  

Let now $D$ be a nodal irreducible cubic; in some coordinates it takes the equation $zy^2=x^3+x^2z$. For every  integer $n>0$, let us define the morphism 
$\Gm\to \Pr_2$ by sending 
$$\Gm\ni t\mapsto (4t^n(t^n-1) : 4t^n(t^n+1) :(t^n-1)^3+8).$$
Again, its image avoids the curve $D$, obtaining the same conclusion as in the cuspidal case.
\medskip

We now come to the most interesting case where $D$ is smooth. We shall produce curves of Euler characteristic $\leq 12$ and of arbitrarily high degree. We let $P_1,\ldots,P_9$ be the flexi of $D$, and $P_{10}, P_{11}, P_{12}$ the three points (apart $P_1$) whose tangent passes through $P_1$. In other words, if $P_1$ is taken as the origin of the elliptic curve corresponding to $D$, then $P_1,\ldots,P_9$ are the points of  three-torsion and $P_{10},P_{11},P_{12}$ are the points of exact order two. 

Put $S=\{P_1,\ldots,P_{12}\}$ and $\c:= D\setminus S$; since $D$ has genus one, $\c$ is an affine curve with $\chi(\c)=12$. 

Finally, let $u\in\O_S^*$ be a non constant $S$-unit; note that $u$ exists and actually $\O_S^*/\kappa^*$ has rank eleven: in fact any point $P_i$, $1<i\leq 12$, is torsion with respect to $P_1$ of order $2$ or $3$,   so there are functions $f_i$ with divisors $(f_i)=6(P_1-P_i)$.  

For every point $p\in D\setminus S$, the tangent $T_p$ at $D$ in $p$ will intersect the cubic $D$ at a second point $p^\prime$, which is distinct from $p$, since $p$ is not a flexus. Also, $T_p$  will intersect the tangent line $T_{P_1}$ at $P_1$ at a point $p^{\prime\prime}$ distinct from $P_1$, since $p\not\in\{P_{10}, P_{11}, P_{12}\}$. The point $p^{\prime\prime}$ does not lie on $D$, since $T_{P_1}\cap D=\{P_1\}$. 

We can canonically identify the line $T_p$ with  $\Pr_1$ in such a way that $p$ corresponds to the point at infinity, $p^\prime$ to $0$ and $p^{\prime\prime}$ to $1$. 

We define a morphism $\c\to\Pr_2$ by sending $p$ to the point $\varphi(p)\in T_p$, where $\varphi(p)$ is the point of coordinate $u(p)$, in the above coordinate system. 

By construction, $\varphi(p)$ never lies on $D$ (otherwise $u(p)=0$ or $\infty$, which does not happen on $\c$). 

The image $\varphi(\c)$ is a curve of degree $\geq\deg(u)$: in fact,  $\varphi(\c)$ intersects $T_{P_1}$ whenever $u(p)=1$ hence the intersection $\varphi(\c)\cap T_{P_1}$ has at least $\deg u$ points, so $\deg \overline{\varphi(\c)}\geq \deg u$. Since $\deg u$ may be taken arbitrarily large (for instance replacing  $u$ with $u^n$) one obtains curves parametrized by $\c$ of arbitrarily high degree. 
\medskip

The result just proved is in accordance with its arithmetic analogue, proved by Beukers in  [Be], stating that the integral points on the complement of a smooth cubic curve on $\Pr_2$ are Zariski dense over a suitable ring of $S$-integers. Its complex-analytic analogue also holds, actually in stronger form: a recent theorem of Buzzard and Lu [Algebraic surfaces holomorphically dominable by $\C^2$, {\it Inv. Math.} {\bf 139} (2000), 617-659] provides the existence of a map $\C^2\to\Pr_2\setminus D$ whose differential has generically rank two; in particular, one obtains Zariski dense  entire curves in the plane omitting the cubic $D$.

\smallskip

Now consider the case when $D$ has two components: \medskip

\noindent{\bf Proposition 4.2}. {\it Let $D\subset\Pr_2$ be the union of a smooth conic
 and a secant line and let $P,Q$ be the singular points of $D$.
 For each positive number $n$ there exists a morphism $f:\Pr_1\rightarrow \Pr_2$
 such that the image $\c:=f(\Pr_1\setminus\{0,\infty\})$ is
contained in $\Pr_2\setminus D$ and the projective curve $f(\Pr_1)$ has degree $n$.
In particular the affine curve $\c\subset\Pr_2\setminus D$ verifies $\deg(\c)=n$
 and $\chi(\c)=0$. }
\medskip

\noindent{\it Proof}. In a suitable system of homogeneous coordinates $(x_0:x_1:x_2)$ ,
the conic will have the equation $x_1^2-x_2^2=x_0^2$ and the line will be defined by $x_0=0$.
Putting $x:=x_1/x_0,\, y=x_2/x_0$, $\ \Pr_2\setminus D$ will be the complement in $\A^2$ of the
hyperbola of equation $x^2-y^2=1$ whose points at infinity $P$ and $Q$ will have coordinates
$P=(0:1:-1)$ and $(0:1:1)$.
Consider, for $n>1$,  the morphism $f:\Gm\rightarrow\A^2$ defined by
$$
f(t)=\left({t^2-t^n+1 \over  2t},{1-t^{n}-t^2\over 2t}\right)=(x(t),y(t))
$$
and observe that $x(t)^2-y(t)^2=1-t^{n}$. Then for no $t\in\Gm$ the point $f(t)$ can
 belong to the curve of equation $x^2-y^2=1$. On the contrary, the continuation of $f$
 to a morphism $\Pr_1\rightarrow\Pr_2$ sends $0$ to the point at infinity $P$ and
$\infty$ to the point at $Q$.
 We have then obtained a curve $\c=f(\Gm)\subset\Pr_2\setminus D$ of zero Euler
characteristic. We now calculate its degree.  Putting for simplification $\xi:=(x+y),\,
\eta:=x-y$ one obtains $\xi=(1-t^{n})/t$ and $\eta=t$, so $\xi$, $\eta$ satisfy the
irreducible equation of degree $n$
$$
\xi\eta=1-\eta^n.
$$
\cvd

When $D$ has three components, we have:
\medskip

\noindent{\bf Proposition 4.3}. {\it Let $D\subset \Pr_2$ be the union of three lines in general position. Then the complement $\Pr_2\setminus D$ is isomorphic to the torus $\Gm^2$. For each positive $n$ there exists a one-dimensional subtorus $\c\subset\Gm^2\simeq \Pr_2\setminus D$,
which has degree $n$ as a curve in $\Pr_2$; its Euler characteristic is $0$.}

\smallskip

\noindent{\it Proof}. Choosing suitable homogeneous coordinates, the divisor $D$ has equation
$x_0x_1x_2=0$. The morphism  $\Gm\ni t\mapsto (1:t:t^n)$ maps $\Gm$ to $\Pr_2\setminus D$. The
image is the degree $n$ curve of equation $x_1^n=x_2$; it is also a one-dimensional
algebraic subgroup with respect to the canonical algebraic group structure of $\Pr_2\setminus
D\simeq  \Gm^2$.
\cvd

\bigskip

We end this paragraph by showing that the condition on the normal crossing singularities of $D$ cannot be removed. In fact we have
\medskip

\noindent{\bf Proposition 4.4}. {\it Let $D$ be the sum of a smooth conic and two lines
 meeting on the conic. Let $X=\Pr_2\setminus D$.
 For every positive $n$ there exists a curve of degree $n$
and vanishing Euler characteristic on $X$.}
\medskip

We can choose affine coordinates $x,y$ for the complement of the first line on $\Pr_2$,
identified with ${\bf A}^2$, so that the second line has equation $x=0$ and the conic
$(x-1)y+1=0$. Then the image of the morphism $f:\Gm\rightarrow \Pr_2\setminus D$ defined by
$$
f(t)=\left(t,{t^{n+1}-1\over t-1}\right)=(x(t),y(t))\
$$
avoids both the line of equation $x=0$ and the conic of equation $(x-1)y+1=0$.
Clearly it is a plane curve of degree $n$.\cvd

 \bigskip

\centerline{\title Appendix}\medskip

The upper bound of  Corollary 2.3 for the number of common zeros
outside $S$ of $a-1,b-1$  is related to the paper [BMZ1],
especially to Thm. 2 therein, which in a different but  equivalent language states:
{\it Let $\c\subset\Gm^n$ be an irreducible curve  over
$\overline\Q$. Suppose that the coordinate functions $x_1,\ldots ,x_n$
on $\c$ are  multiplicatively independent modulo constants.
Then there are only finitely many points $P\in \cc$ such that there are
two independent relations $x_1(P)^{a_1}\cdots x_n(P)^{a_n}=x_1(P)^{b_1}\cdots x_n(P)^{b_n}=1$,
$a_i,b_i\in\Z$.} (See [BMZ2] for the case over $\bf C$.)

On setting   $a=x_1^{a_1}\cdots x_n^{a_n}$, $b=x_1^{b_1}\cdots x_n^{b_n}$,
we see that each relevant point $P$ in this statement is a common zero of $a-1,b-1$,
clarifying the connection with the present context.

For instance Corollary 2.3 (for $\kappa=\overline{\Q}$, $S=$ set
of zeros/poles of the $x_i$) bounds the number of relevant points $P$, for a
{\it given} choice of the exponents $a_i,b_i$. In a way,  [BMZ1, Thm. 2] is
much stronger, in that it asserts the finiteness of the points  for {\it varying}
 pairs of exponents vectors,   leading for instance to a bound
$\sum_{v\not\in S}\min (v(1-a),v(1-b))=O_{\cc,S}(1)$, independently of $a,b$ of the above shape.
However the present corollary is rather more explicit and uniform in the dependence
with respect to $\cc$ and $S$.
This may be crucial for some applications; for instance, the result of [BMZ1]
 would not be sufficient to derive the present Theorems 1.1 and 1.3;
in fact, in the present proofs Corollary 2.3 is applied with a set $S$
(and also a curve $\c$) which  vary with the individual solutions,
and  hence a good uniformity is needed.\smallskip

Also, Corollary 2.3 would lead to a simplification of the proofs in [BMZ1]; we
briefly indicate how, referring to that paper for this argument. We let $a,b$ be as above.
Then  Corollary  2.3 gives
$\sum_{v\not\in S}\min(v(1-a),v(1-b))\ll (H(a)H(b))^{1/3}$.

Say now  that $\c$ is defined over the number field
$L$. Then if $a(P)-1=b(P)-1=0$ for a $P\in\c$, the same is true for any conjugate of $P$
over $L$. Then  Corollary 2.3 implies at once $[L(P):L]\ll (H(a)H(b))^{1/3}$.
We may also assume that the exponent vectors $\a,\b$ are the first
 two successive minima for the lattice of exponents in the  relations among the $x_i(P)$,
 with volume   $V$, say. We then  have $H(a)H(b)\ll|\a|\cdot|\b|\ll V$. Hence $[L(P):L]\ll
V^{1/3}$. Now, in the notation of [BMZ1], we may take $V\ll N\Pi$,  and then the
above  improves on [BMZ1, (4.4)], in the crucial case $r=n-2$; the gain suffices   to avoid the
final, somewhat involved,  arguments in that proof.\medskip

 We should mention here also the result by Ailon and Rudnick [AR]:
they consider the $\gcd(f^n-1,g^n-1)$ for given polynomials $f,g$
(or more generally for given rational functions on a curve).
Using the well-known finiteness of torsion points on curves in $\Gm^n$ not contained
 in any translate  of a proper algebraic subgroups,
they show that the $\gcd$ is ``usually" $1$, for multiplicatively independent $f,g$.
\medskip
As in the previous situation concerning [BMZ1], a very special case of Corollary 2.3
 bounds the $\gcd$-degree for given
$n$, which is on the one hand weaker than [AR].
 However, somewhat surprisingly, this bound allows to recover easily the arithmetical result
 about torsion points. In fact, let
$x,y$ be multiplicatively independent rational functions on the curve
$\c/L$ and let $P$ be a torsion point, so $x(P)^N=y(P)^N=1$ for some integer $N\ge 1$
 which we suppose to be minimal. Then every conjugate of $P$ over $L$ is a common zero
of $x^N-1 $ and $y^N-1$. Applying Corollary 2.3 we deduce that $[L(P):L]\ll N^{2/3}$.
But since $N$ is minimal, the relevant degree is $\gg \phi(N)\gg N/\log\log N$,
 proving the boundedness of $N$.

The argument actually leads to quite explicit bounds (in terms
only of the degree and genus of $\c$ and of $[L:\Q]$) whose
calculation we shall perform in a future note. Here we only remark
that they lead to:

 (i) recover the best possible exponent (i.e. $2$) for  the
degree and

(ii) better results than the known ones when the genus is small.

To our knowledge no known estimate involves a dependence on the
genus.
\medskip

It seems not a very common fact  that  function-field arguments like the present ones,
 involving derivations, lead directly to
 arithmetical deductions; it would be interesting to develop function field
 techniques for similar questions concerning abelian varieties;
apart from some independent interest, this could perhaps lead to analogous, much deeper,
arithmetical consequences (first obtained in this direction by Raynaud).
\bigskip

{\bf Acknowledgements}. The authors are grateful to Professors Michel Raynaud and Stephen Yau and to an anonimous referee for their consideration of the present work and several suggestions leading to a better  presentation of the paper.
\bigskip

$${}$$

{\bf References}.
\medskip

\item{[AR]} - N. Ailon, Z. Rudnick, Torsion points on curves and common divisors of
$a^k-1$ and  $b^k-1$, {\it Acta Arith.} {\bf 113} (2004), 31-38.  \smallskip

\item{[Be]} - F. Beukers, Ternary Form Equations, {\it J. Number Theory} {\bf 54} (1995), 113-133.
\smallskip

\item{[BMZ1]} - E. Bombieri, D. Masser, U.Zannier, Intersecting a curve with algebraic
subgroups of multiplicative groups, {\it Int. Math. Research Notices} {\bf 20} (1999),
1119-1140.\smallskip

\item{[BMZ2]} - E. Bombieri, D. Masser, U.Zannier, Finiteness results for multiplicative
dependent points on complex curves, {\it Michigan Math. J.} {\bf 51} (2003), 451-466.\smallskip

\item{[Bo]} - E. Borel, Sur les z\'eros des fonctions enti\`eres, {\it Acta Math.}, {\bf 20}
(1897), 357-396.\smallskip

\item{[BrM]} - D. Brownawell, D. Masser, Vanishing sums in function fields,
{\it Math. Proceedings Cambridge Phyl. Soc}, {\bf 100} (1986), 427-434.\smallskip

\item{[CoZ]} - P.B. Cohen, U. Zannier, Fewnomials and Intersections
of Lines with Real Analytic Subgroups in $\Gm^n$,
{\it Bull. London Math. Soc.} {\bf 34} (2002), 21-32.\smallskip

\item{[CZ1]} - P. Corvaja, U. Zannier, On the diophantine equation $f(a^m,y)=b^n$,
{\it Acta Arithmetica} {\bf 94.1} (2000), 25-40.\smallskip

\item{[CZ2]} - P. Corvaja, U. Zannier,
A lower bound for the height of a rational function at
$S$-unit points, {\it Monatshefte f. Math.} {\bf 144} (2005), 203-224.
\smallskip

\item{[G1]} - M. Green, On the functional equation $f^2=e^{2\phi_1}+
e^{2\phi_2}+e^{2\phi_3}$ and a new Picard theorem, {\it Transactions
of the American Math. Soc.}, {\bf 195} (1974), 223-230.\smallskip

\item{[G2]} - M. Green, Some Picard theorems for holomorphic maps
to algebraic varieties, {\it American J. Math.}, {\bf 97} (1975), 43-75.
\smallskip

\item{[KMK]} - S. Keel, J. McKernan, Rational Curves on Quasi Projective
Surfaces, {\it Memoirs of the American Math. Soc.}, {\bf 669}, (1999).\smallskip

\item{[NWY]} - J. Noguchi, J. Winkelmann, K. Yamanoi, The second main theorem for  holomorphic curves into semi-abelian varieties II, Forum Math.{\bf 20} (2008),  469--503.
\smallskip

\item{[S]} - J.H. Silverman, Generalized Greatest Common
Divisor, Divisibility Sequences, and Vojta's Conjecture for Blowups,
{\it Monatsh. f. Math.} {\bf 145} (2005), 333-350.\smallskip

\item{[V]} - P. Vojta, Diophantine Approximations and Value Distribution Theory, LNM 1239, Springer
1987.\smallskip

\item{[W]} - J.T. Wang, An effective Schmidt's subspace theorem over function fields,
{\it Math. Zeit.}, {\bf 246} (2004), 811-844.\smallskip

\item{[Z1]} - U. Zannier, Some remarks on the $S$-unit equation in function
fields, {\it Acta Arith.} {\bf 94} (1993), 87-98.\smallskip

\item{[Z2]} - U. Zannier, Polynomial squares of the form
$aX^m+b(1-X)^n+c$, {\it Rend. Sem. Mat. Univ. Padova}, {\bf 112}
(2004), 1-9.\smallskip

\bigskip

$${}$$

Pietro Corvaja\hfill Umberto Zannier

Dipartimento di Matematica e Informatica\hfill  Scuola Normale Superiore

Via delle Scienze, 206\hfill Piazza dei Cavalieri, 7

33100 - Udine (ITALY)\hfill  56100 Pisa (ITALY)

corvaja@dimi.uniud.it\hfill u.zannier@sns.it

\end

\bye